\documentclass{elsarticle}
\usepackage[english]{babel}
\usepackage[utf8x]{inputenc}
\usepackage[T1]{fontenc}
\usepackage{graphicx}
\usepackage[colorinlistoftodos]{todonotes}
\usepackage[colorlinks=true, allcolors=tudelftblue]{hyperref} 
\usepackage{amssymb}
\usepackage{caption}
\usepackage{subcaption}
\usepackage{xcolor}
\usepackage{float}
\usepackage{algorithm}
\usepackage[noend]{algpseudocode}
\usepackage{titling} 
\usepackage{blindtext}\usepackage{titlesec}
\usepackage[colorinlistoftodos]{todonotes}
\usepackage{tikz}
\usepackage{array}
\usepackage{geometry}
\usepackage{mathtools}
\usepackage{hyperref} 
\usepackage{sectsty}
\usepackage{amsmath}
\usepackage{listings}  
\usepackage{color}     
\usepackage{pgfplots}
\pgfplotsset{compat=newest}
\usepackage{amsthm}

\newtheorem{theorem}{Theorem}
\newtheorem{assumption}{Assumption}
\newtheorem{remark}{Remark}
\usepackage{amsmath}
\usepackage{tikzpagenodes}
\usepackage[justification=centering]{caption}
\usepackage{booktabs}
\usepackage{listings}
\numberwithin{equation}{section}
\definecolor{tudelftdarkblue}{RGB}{12,35,64}
\definecolor{tudelftcyan}{RGB}{0,166,214}
\definecolor{tudelftblue}{RGB}{0, 118, 194}
\allsectionsfont{\color{black}} 
\usepackage{accents}
\newcommand{\ubar}[1]{\underaccent{\bar}{#1}}
\newcommand{\customfootnotetext}[2]{{
  \renewcommand{\thefootnote}{#1}
  \footnotetext[0]{#2}}}
\titleformat{\paragraph}
{\normalfont\normalsize\bfseries}{\theparagraph}{1em}{}
\titlespacing*{\paragraph}
{0pt}{3.25ex plus 1ex minus .2ex}{1.5ex plus .2ex}
\usepackage{setspace}
\begin{document}

\newpage
\begin{center}\customfootnotetext{$\dagger$}{Electronic Supplementary Information (ESI) available}
{\LARGE INTHOP: A Second-Order Globally Convergent Method for Nonconvex Optimization\textsuperscript{$\dagger$}}
\end{center}

\begin{center}
{\textit{Krishan Kumar\textsuperscript{1,\#} Ashutosh Sharma\textsuperscript{1,\#}\customfootnotetext{\#}{These authors contributed equally}, Gauransh Dingwani\textsuperscript{2}, Nikhil Gupta\textsuperscript{1}, Vaishnavi Gupta\textsuperscript{1} and Ishan Bajaj\textsuperscript{1,*}\customfootnotetext{*}{Corresponding author: email: \href{mailto:ibajaj@iitk.ac.in}{ibajaj@iitk.ac.in}}}}
\\
{\small \textit{$^1$Department of Chemical Engineering, Indian Institute of Technology Kanpur, Kanpur, Uttar Pradesh, India}}\\
{\small \textit{$^2$Department of Chemical Engineering, Indian Institute of Technology Roorkee, Roorkee, Uttarakhand, India}}
\addcontentsline{toc}{section}{Abstract}
\section*{Abstract}
\begin{quote}
Second-order Newton-type algorithms that leverage the exact Hessian or its approximation are central to solve nonlinear optimization problems. These algorithms have been proven to achieve a faster convergence rate than the first-order methods and can find second-order stationary points. However, their applications in solving large-scale nonconvex problems are hindered by three primary challenges: (1) the high computational cost associated with Hessian evaluations, (2) its inversion, and (3) ensuring descent direction at points where the Hessian becomes indefinite. We propose INTHOP, an interval Hessian-based optimization algorithm for nonconvex problems to deal with these primary challenges. We propose a new search direction by approximating the original Hessian by a positive definite matrix valid within an interval. This results in always a descent direction and unlike other second-order methods, our method does not require Hessian computation at each iteration. We prove that the difference between the positive definite approximate and exact Hessian is bounded within an interval. We develop various algorithm variants based on the interval size updating methods and Hessian approximation computation methods. We also prove the global convergence of the proposed algorithm. {Further, we apply the algorithm to an extensive set of test problems, including high dimensional cases, and compare its performance with the existing methods such as steepest descent, limited memory Broyden–Fletcher–Goldfarb–Shanno (LBFGS) algorithm implemented in PyTorch, and the open-source solver IPOPT}. {We show empirically that INTHOP demonstrates competitive performance with IPOPT and significantly outperforms LBFGS and steepest descent across all metrics}.
\end{quote}
\end{center}
\section{Introduction}
We consider the following unconstrained optimization problem: 
\begin{equation} \label{eqn:intro}
    \min_{\mathbf{x} \in \mathbb{R}^{n}} f(\mathbf{x})
\end{equation}
where \( f: \mathbb{R}^n \to \mathbb{R} \) is a twice continuously differentiable but possibly nonconvex function and $n$ is the number of variables. Over the past decades, several zero-, first- and second-order methods based on using gradient and Hessian information, respectively, have been developed for problem \eqref{eqn:intro}, for instance, see {\cite{cartis2010complexity, gratton2025yet, liu1989limited, nesterov2006cubic, Shanno1970ConditioningOQ, zhang1999new}} and references therein. First-order methods are often preferred to solve large-scale optimization problems in various areas, including engineering, machine learning, image processing, and computational chemistry. Steepest descent is one of the most popular first-order algorithms for generating a sequence of iterates $\mathbf{x}_k$ where $k \in \{0, 1, 2, \dots, p\}$ such that \(||\nabla f(\mathbf{x}_p)|| \xrightarrow{} 0\). Due to the slow convergence rate of the steepest descent algorithm, a class of effective first-order algorithms has been developed, which has drawn the attention of the optimization community. Some of these first-order methods are batch gradient descent, stochastic gradient descent \cite{robbins1951stochastic}, mini-batch gradient descent, Momentum \cite{qian1999momentum,liu2020improved}, Nesterov accelerated gradient \cite{ghadimi2016accelerated,1370862715914709505}, Adagrad, RMSprop, and Adam \cite{kingma2014adam, mcmahan2010adaptive, zeiler2012adadeltaadaptivelearningrate}.

However, first-order algorithms have several drawbacks.  They converge slowly for high-dimensional problems. Also, they do not adequately capture the topology of the landscape. 
To this end, second-order methods are promising and have shown, both empirically and theoretically, a faster convergence rate to second-order stationary points and resilience to ill-conditioning of the problem. The search direction in a second-order method is based on the Hessian or its approximation and gradient information of the objective function. Computing the Hessian enables capturing the function's curvature. Newton's method is a popular second-order method to solve nonlinear optimization problems. 

Although Newton's method exhibits impressive convergence characteristics compared to first-order approaches \cite{ karimireddy2018global, nocedal1999numerical}, its broad application to large-scale nonconvex optimization problems poses significant challenges. First, the memory requirement increases with the square of the number of variables. Second, the Hessian calculation can be expensive. Third, finding the search direction requires Hessian matrix inversion at each iteration. Finally, Newton's direction for nonconvex problems may not be descent at each iteration. However, several practical applications involve nonconvex problems, including training neural networks, computational chemistry, designing distillation columns, and reactor networks. 

Four broad classes of second-order algorithms have been developed to overcome one or more challenges associated with the Newton method: quasi-Newton, sub-sampled Newton, regularized Newton, and Newton method with lazy Hessian updates. Quasi-Newton methods, such as Broyden-Fletcher-Goldfarb-Shanno (BFGS) method, see \cite{Broyden1970TheCO,Fletcher1970ANA,Goldfarb1970AFO,Shanno1970ConditioningOQ}, overcome the challenges of Hessian computation by approximating the Hessian using gradients of the previous iterations while ensuring the approximate Hessian matrix is positive definite. Sub-sampled Newton methods estimate Hessian using a subset of data to solve finite sum minimization problems that commonly occur in training machine learning models, for instance, see \cite{bollapragada2016exactinexactsubsamplednewton, erdogdu2015convergence,roosta2019sub}. A class of second-order algorithms based on cubic regularization \cite{nesterov2006cubic} and regularizing the Hessian matrix with the square root of the gradient has also been proposed \cite{gratton2025yet} and shown to have attractive convergence rates \cite{cartis2010complexity}. The regularization parameter is selected to ensure that the matrix is positive definite and a certain reduction is obtained at each iteration. Newton method with lazy Hessian updates avoids Hessian computation for several iterations while computing gradients at each iteration \cite{doikov2023second}. 
Several contributions have been made based on exploiting the problem structure \cite{allman2019decode, mitrai2021efficient} and using decomposition approaches \cite{chiang2014structured, kang2014interior, li2019generalized} and GPU architectures \cite{pacaud2024parallel, montoison2025madncl} to solve large-scale nonlinear programming problems. 

The BFGS method has the advantages that the Hessian computation is avoided and the approximate Hessian is positive definite, thus ensuring that the search direction is always descent. The method has a drawback in that it requires satisfying the curvature condition, which is satisfied only for convex problems. For nonconvex problems, satisfying the curvature condition requires finding a step length that satisfies computationally expensive strong Wolfe conditions. Sub-sampled Newton methods require less storage, but they are valid for only convex problems involving finite sum minimization. Regularized Newton methods have been shown to converge to a first-order minimizer in at most $\mathcal{O}(|log\ \epsilon|\epsilon^{-3/2})$ function evaluations. The disadvantage of the method is that the eigenvalues and Lipschitz constant for the Hessian need to be computed at each iteration. Newton method with lazy Hessian updates has been shown to improve on the complexity of cubic Newton by a factor of $\sqrt{n}$. A major drawback of the method is that while a simple Newton step requires solving a linear equation, more complex procedures, such as Lanczos pre-processing \cite{cartis2011adaptive}, are required to find the next iterate.   

{ In summary, despite significant progress in second-order optimization algorithms, existing approaches typically address only one of the following challenges at a time: (i) avoiding full Hessian computation (quasi-Newton), (ii) enforcing positive definiteness through per-iteration regularization (regularized Newton), or (iii) skipping Hessian updates according to a predefined schedule (lazy Newton). \textit{In contrast, we propose a fundamentally different approach based on the Hessian approximation in an interval}. Specifically, instead of modifying or approximating the Hessian at each iteration, we construct a positive definite Hessian approximation that is valid over an entire interval.  This is achieved by calculating the Hessian of the $\alpha$BB convex underestimator \cite{adjiman1998global2} of the objective function at any point inside the interval region, i.e., $[\mathbf{x}^L, \mathbf{x}^U]$, where $
[\mathbf{x}^L, \mathbf{x}^U]
=
\left\{
\mathbf{x} \in \mathbb{R}^n : 
x_i^L \le x_i \le x_i^U,\; \text{ for all } i = 1,2, \ldots, n
\right\}
$. Therefore, the main novelty of the proposed method lies in constructing a region-wise positive definite Hessian approximation that guarantees the proposed search direction to be descent while requiring Hessian evaluation and inversion only at specific iterations. \textit{To the best of our knowledge, this is the first approach that constructs a positive definite Hessian approximation valid over an interval region in the context of nonconvex optimization}. 
}  
Our analysis shows that the difference in the Hessian of the $\alpha$BB underestimator at a point $\mathbf{z}\in[\mathbf{x}^L, \mathbf{x}^U]$ and Hessian of the original function for all $\mathbf{x} \in [\mathbf{x}^L, \mathbf{x}^U]$ ($\nabla^2f(\mathbf{x})$) is bounded and proportional to the interval size. We incorporate the search direction in a line-search framework and develop an algorithm, INTHOP - INTerval Hessian-based OPtimization. We develop two variants of the INTHOP algorithm based on fixed and adaptive interval sizes. Our algorithm is similar in spirit to the Newton method with lazy Hessian updates in that we do not update the Hessian at each iteration. However, our algorithm is different in two aspects. Firstly, instead of updating the Hessian matrix after a fixed schedule, we update it based on the interval for which the difference between the approximate Hessian and the exact Hessian is bounded. Secondly, finding the search direction is computationally inexpensive and can be obtained by taking a matrix-vector product or solving a linear equation. We prove the global convergence of INTHOP in the sense that there exists an accumulation point of the sequence generated by INTHOP that is a stationary point. Through extensive computational experiments, we illustrate the effectiveness of various variants of INTHOP. We also compare INTHOP with { IPOPT}, steepest descent, and  LBFGS. 

The rest of the article is organized as follows. Section 2 introduces the notations. Section 3 outlines the proposed methodology, starting with the construction of the search direction, followed by the formulation and illustration of the interval Hessian. The overall algorithm structure and the adaptive delta strategy are also detailed here. Section 4 briefly describes the competing algorithms used for comparison. Section 5 presents numerical experiments and results, beginning with illustrative examples and proceeding to a comprehensive benchmarking of the interval Hessian frameworks under different conditions, concluding with a comparison against other algorithms. Finally, Section 6 provides the conclusion of the study.

\section{Notations}
Throughout this article, a vector of size $n$ is denoted using bold lowercase letters, e.g., $\mathbf{a} = (a_1, a_2, \dots, a_n)^\top$ where $a_i$ is the $i^{\text{th}}$ component of the vector $\mathbf{a} \in \mathbb{R}^n$ and $\mathbf{a}^\top$ denotes its transpose. $\mathbf{x}_k$ is the iterate generated at iteration $k$. For matrices, we use uppercase bold letters, e.g., $\mathbf{A} = (a_{ij})$ where $a_{ij}$ is the component in $i^{\text{th}}$ row and $j^{\text{th}}$ column. A matrix is symmetric when $a_{ij} = a_{ji}$. The notation $\|.\|$ represents the $\ell_2$ norm. The $i^{\text{th}}$ eigenvalue of \(\mathbf{A}\) is denoted by \(\lambda_i(\mathbf{A})\). The symbol $\rho(\mathbf{A})$ represents the spectral radius of matrix $\mathbf{A}$. We use interval variables denoted by lowercase letters in brackets, e.g. $[a]$, with the corresponding lower and upper bounds denoted by $\underline{a}$ and $\overline{a}$, respectively. The width of the interval variable is given by $w([a]) = \overline{a} - \underline{a}$. An interval vector is given by a bold lowercase letter enclosed in square brackets, e.g., $[\mathbf{a}] = \left([a_1], [a_2], \dots, [a_n]\right)^\top, [a_i] = [\underline{a}_i, \overline{a_i}]$. The width of the interval vector $[\mathbf{a}]$ is given by $w([\mathbf{a}]) = \max_{i} (\overline{a}_i - \underline{a}_i).$ Similarly, we represent an interval matrix using bold uppercase letters enclosed in square brackets, e.g., $[\mathbf{A}] = ([{a}_{ij}]), {a}_{ij} = [\underline{a}_{ij}, \overline{a}_{ij}]$. An interval matrix is symmetric when $[{a}_{ij}] = [{a}_{ji}]$ \cite{moore2009introduction}.

Let $f: \mathbb{R}^n \to \mathbb{R}$ be a twice continuously differentiable function. We denote the gradient of $f$ by $\nabla f(\mathbf{x}) \in \mathbb{R}^n$ and the Hessian of $f$ by $\nabla^2 f(\mathbf{x}) \in \mathbb{R}^{n \times n}$. We use the following notation throughout: $f_k:= f(\mathbf{x}_k), \mathbf{g}_k:= \nabla f(\mathbf{x}_k):= \nabla f_k, \text{ and } \nabla^{2} f_{k} := \nabla^2 f(\mathbf{x}_k)$. 

\section{Methodology}
\subsection{Search Direction Construction}
\label{sec:SearchDirectionConstruction}
In this article, we propose a variant of Newton direction that ensures descent for nonconvex problems at each iteration. Classical Newton's method is based on approximating $f(\mathbf{x}_k+\mathbf{p})$ by a second-order Taylor series model as follows:
\begin{equation}
    f(\mathbf{x}_k+\mathbf{p}) \approx m_k^N(\mathbf{p}) = f_k + \nabla f^{\top}_k\mathbf{p} + \frac{1}{2}\mathbf{p}^{\top}\nabla^2f_k \mathbf{p}.
\end{equation}
The Newton direction is obtained by minimizing $m_k^N(\mathbf{p})$ and  is given by
\begin{equation}
    \mathbf{p}_k^N = -(\nabla^2 f_k)^{-1}\nabla f_k \text{ whenever } (\nabla^2 f_k)^{-1} \text{ exists}.
\end{equation}
If matrix $\nabla^2 f_k$ is not positive definite, then $\mathbf{p}_k^N$ may not satisfy the descent property that is $\nabla f_k^{\top}\mathbf{p}_k^N < 0$. Accordingly, we propose the following quadratic model to approximate $f(\mathbf{x}_k+\mathbf{p})$:
\begin{equation}
    f(\mathbf{x}_k+\mathbf{p}) \approx m_k^{IH}(\mathbf{p}) = f_k + \nabla f^{\top}_k\mathbf{p} + \frac{1}{2}\mathbf{p}^{\top}\nabla^2\mathcal{L}_t \mathbf{p}
    \label{eq:Modified-Model}
\end{equation}
\begin{equation}
    \mathbf{p}_k^{IH} = -(\nabla^2\mathcal{L}_t)^{-1}\nabla f_k \text{ whenever } (\nabla^2 \mathcal{L}_{t})^{-1} \text{ exists, } 
    \label{eqn:searchdirection}
\end{equation}
where
\begin{equation}
    \mathcal{L}(\mathbf{x}) = f(\mathbf{x}) + \sum_{i=1}^n \alpha_i({x}^L_i - {x}_i)({x}^U_i - {x}_i).
    \label{fig:UnderEstimator}
\end{equation}
\begin{equation*}
    \mathbf{x}, \mathbf{x}_t,\mathbf{x}_k \in [\mathbf{ x}^L,\mathbf{x}^U] \text{ and } \alpha_i \geq 0.
\end{equation*}
Here, $\mathcal{L}(\mathbf{x})$ is the $\alpha$BB underestimator used to solve nonconvex problems to global optimality \cite{adjiman1998global2,androulakis1995alphabb}. In the above underestimator, a negative term is added to $f(\mathbf{x})$, ensuring $\mathcal{L}(\mathbf{x})\leq f(\mathbf{x}), \forall \mathbf{x} \in [\mathbf{x}^L,\mathbf{x}^U]$. The $\alpha_i$ for each $i=1,2,\ldots,n$ is computed such that $\mathcal{L}(\mathbf{x})$ is guaranteed to be convex for each $\mathbf{x} \in [{\mathbf{x}^L},\mathbf{x}^U]$. This implies that the Hessian of $\mathcal{L}(\mathbf{x})$ is positive semi-definite. The Hessian of $\mathcal{L}(\mathbf{x})$ and $f(\mathbf{x})$ is related as follows
\begin{equation}
    \nabla^2\mathcal{L}(\mathbf{x}) = \nabla^2f(\mathbf{x}) + 2\Delta,
    \label{fig:undrest}
\end{equation}
where $\Delta$ is a diagonal matrix of order $n$ whose diagonal elements are $\alpha_i$ for all $i=1,2,\ldots,n$. The matrix $\Delta$ is referred as the diagonal shift matrix. For simplicity, we assume a uniform diagonal shift matrix, i.e. $\alpha_1=\alpha_2=\dots=\alpha_n=\alpha$, which makes the above relation as follows:
\begin{equation}
    \nabla^2\mathcal{L}(\mathbf{x}) = \nabla^2f(\mathbf{x}) + 2\alpha {I},
    \label{eq:commonalpha}    
\end{equation}
where ${I}$ is an identity matrix of order $n$. We denote an identity matrix of order $n$ by $I$ throughout this work. Maranas and Floudas \cite{androulakis1995alphabb} proved that the convex underestimator $\mathcal{L}(\mathbf{x})$, defined in \eqref{fig:UnderEstimator}, is convex if and only if
\begin{equation}
    \alpha \geq \max\left\{{0, -\frac{1}{2}\min_{i,\mathbf{x}^L\leq \mathbf{x} \leq \mathbf{x}^U}\lambda_i(\nabla^2f(\mathbf{x}))}\right\},
    \label{eq:alpha}
\end{equation}
where  $\lambda_i$'s are the eigenvalues of $\nabla^2f(\mathbf{x})$. This information will be a key factor in this work.

Next, we present a result that shows that the differences between $L_t$ and $f_k$, $\nabla L_t$ and $\nabla f_k$, and $\nabla^2L_t$ and $\nabla^2 f_k$ are bounded. 
\begin{theorem}\label{theorem:theo1}
Let $f: \mathbb{R}^n \rightarrow \mathbb{R}$ be a twice continuously differentiable function, and $\mathcal{L}(\mathbf{x})$ is an $\alpha$BB convex underestimator constructed on an interval of width $\delta$ such that $\mathbf{x}_t, \mathbf{x}_k\in [\mathbf{x}^L, \mathbf{x}^U]$. Suppose $f$, $\nabla f$, and $\nabla^2 f$ are Lipschitz continuous over a domain $\Omega\subseteq \mathbb{R}^n$. That is, $\forall \mathbf{x}_k, \mathbf{x}_t \in \Omega$, there exist constants $L_f$, $L_g$, and $L_H$ such that
\begin{equation}\label{eq:lipsch_f}
    \|f(\mathbf{x}_t) - f(\mathbf{x}_k)\| \leq L_f \|\mathbf{x}_t - \mathbf{x}_k\|    
\end{equation}
\begin{equation}\label{eq:lipsch_g}
    \|\nabla f(\mathbf{x}_t) - \nabla f(\mathbf{x}_k)\| \leq L_g \|\mathbf{x}_t - \mathbf{x}_k\|    
\end{equation}
\begin{equation}\label{eq:lipsch_h}
    \|\nabla^2 f(\mathbf{x}_t) - \nabla^2 f(\mathbf{x}_k)\| \leq L_H \|\mathbf{x}_t - \mathbf{x}_k\|.
\end{equation}
Suppose the minimum eigenvalue of $\nabla^2 f(\mathbf{x})$ for $\mathbf{x}\in[\mathbf{x}^L, \mathbf{x}^U]$ is $\lambda_{min}$. Then, the following bounds hold:
\begin{align}
    \left| \mathcal{L}(\mathbf{x}_t) - f(\mathbf{x}_k) \right| & \leq \frac{L}{2} \sqrt{n} \delta + \frac{|\lambda_{min}|}{8} n \delta^2, \\
    \left\| \nabla \mathcal{L}(\mathbf{x}_t) - \nabla f(\mathbf{x}_k) \right\| &\leq \frac{L_g}{2} \sqrt{n} \delta + \frac{|\lambda_{min}|}{2} \sqrt{n}\delta, \\
    \left\| \nabla^2 \mathcal{L}(\mathbf{x}_t) - \nabla^2 f(\mathbf{x}_k) \right\| &\leq \frac{L_H}{2} \sqrt{n} \delta + \frac{|\lambda_{\min}|}{2} \sqrt{n}\label{boundlx}. 
\end{align}
\end{theorem}
\begin{proof}
Refer ~\ref{appendix:proofs}. 
\end{proof}
 
{T}wo key features of the $\alpha$BB underestimator that enable us to employ its Hessian in $m_k^{IH}$ are that it is twice differentiable and convex, which ensures that $\nabla^2\mathcal{L}_t$ is guaranteed to be positive semidefinite. Since the Hessian matrix must be positive definite to ensure that the direct is descent,  therefore to make $\nabla^2\mathcal{L}_t$ positive definite, we add a positive term $c_1\tilde{g}$ in it. We denote this new matrix by $\nabla^{2}\mathcal{H}(\mathbf{x})$, i.e, $\nabla^{2}\mathcal{H}(\mathbf{x})=\nabla^{2}\mathcal{L}(\mathbf{x})+c_1\tilde{g}I$, where $c_1$ and $\tilde{g}$ are real positive numbers.

\begin{remark}
It can be easily seen that the result proved in \eqref{eq:lipsch_h} in Theorem 1 is also applicable for $\nabla^{2}\mathcal{H}(\mathbf{x})$ with a different bound. Specifically, in this case, we have 
\begin{equation}\label{boundhx}
\left\| \nabla^2 \mathcal{H}(\mathbf{x}_t) - \nabla^2 f(\mathbf{x}_k) \right\| \leq \frac{L_H}{2} \sqrt{n} \delta + \frac{|\lambda_{\min}|}{2} \sqrt{n}+c_1\tilde{g}\sqrt{n}.
\end{equation} 
\end{remark}

We consider the following univariate nonconvex function to illustrate the behavior of the quadratic models ($m_k^N$ and $m_k^{IH}$) at a point where the Hessian is indefinite. 
\[
f(x) = x^4 - 3x^3 - 1.5x^2 + 10x.
\]

\begin{figure}[H]
\centering
\includegraphics[width=0.5\textwidth]{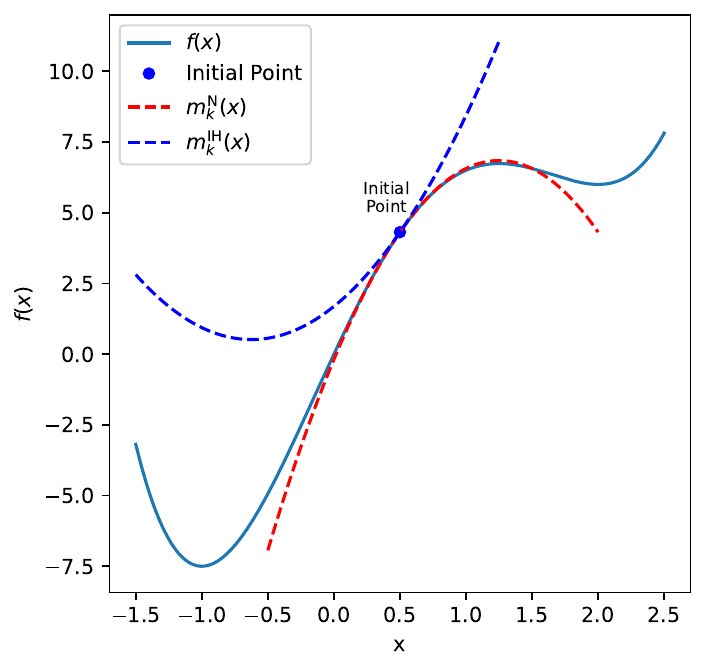}
\caption{Graphical illustration of $f(\mathbf{x}),\ m_k^{N}(\mathbf{x})\ \text{and} \ m_k^{IH}(\mathbf{x})$}
\label{fig:quadModel}
\end{figure}

In Figure~\ref{fig:quadModel}, we compare the quadratic model \( m_k^N({x}) \), constructed using the second-order Taylor expansion of \( f({x}) \), with a modified model \( m_k^{IH}({x}) \) that uses the Hessian of $\mathcal{H}(x)$. We observe that $p_k^N$ would result in increasing the objective function and $p_k^{IH}$ would lead to a descent direction. This example illustrates the importance of a positive definite matrix in ensuring a descent direction. 
Computing $\alpha$ exactly over an interval using Eq. \ref{eq:alpha} is equivalent to solving a nonconvex optimization problem to global optimality, which could be more difficult than solving Eq. \ref{eqn:intro}. Therefore, an efficient and scalable alternative is needed to obtain a reliable lower bound on the minimum eigenvalue over the interval. In the following subsection, we introduce the concept of the interval Hessian \cite{moore2009introduction} and explain how it can be leveraged to compute a bound on the minimum eigenvalue over an interval in a computationally tractable manner. 

\subsection{Interval Hessian} \label{sec: int_Hessian}
We adopt an interval-based Hessian approximation strategy to overcome the computational challenges associated with estimating $\alpha$ using Eq. \ref{eq:alpha}. Specifically, we replace the exact variable-dependent Hessian ($\nabla^2\mathcal{H}(\mathbf{x})$) with an interval Hessian matrix ($[\nabla^2\mathcal{H}(\mathbf{x})]$) computed using interval arithmetic such that $\nabla^2\mathcal{H}(\mathbf{x}) \subseteq [\nabla^2\mathcal{H}(\mathbf{x})]$  for all $\mathbf{x}\in [\mathbf{x}^L, \mathbf{x}^U]$. Such strategy offers considerable scalability and computational advantages.

This interval Hessian strategy introduces a fundamental trade-off: the accuracy of the Hessian approximation versus computational cost. Specifically, using intervals naturally leads to conservativeness—looser intervals yield less accurate Hessian approximation (Theorem 1) and thus potentially less effective search directions. Conversely, tighter intervals reduce this conservativeness, but require increased computational effort because the Hessian needs to be recomputed and the linear equation (\ref{eqn:searchdirection}) resolved to find the search direction if the iterate $\mathbf{x}_k \notin [\mathbf{x}^L,\mathbf{x}^U]$ because $||\nabla^2 \mathcal{L}(\mathbf{x}_{k-1}) - \nabla^2 f(\mathbf{x}_{k})||$ may not be bounded. 

\subsubsection{Interval Hessian Illustrative example}
We illustrate interval Hessian computation using the following problem:
\begin{equation}
    f({x}_1,{x}_2) = \left(1.5 - {x}_1(1 - {x}_2)\right)^2 + \left(2.25 - {x}_1(1 - {x}_2^2)\right)^2 + \left(2.625 - {x}_1(1 - {x}_2^3)\right)^2    
    \label{eqn:example}
\end{equation}
\[
0 \leq {x}_1, {x}_2 \leq 2.
\]
Then,
\[
\begin{aligned}
\frac{\partial^2 f}{\partial {x}_1^2} &= 2(1 - {x}_2)^2 + 2(1 - {x}_2^2)^2 + 2(1 - {x}_2^3)^2 \\
\frac{\partial^2 f}{\partial {x}_2^2} &= {x}_1\left(30{x}_1 {x}_2^4 + 12{x}_1 {x}_2^2 - 12{x}_1 {x}_2 - 2{x}_1 + 31.5{x}_2 + 9\right) \\
\frac{\partial^2 f}{\partial {x}_1 \partial {x}_2} &= 12{x}_1 {x}_2^5 + 8{x}_1 {x}_2^3 - 12{x}_1 {x}_2^2 - 4{x}_1 {x}_2 - 4{x}_1 + 15.75{x}_2^2 + 9{x}_2 + 3.
\end{aligned}
\]
We perform interval arithmetic using the Gaol package \cite{goualard2005gaol} to obtain the following interval Hessian.
\begin{equation}
    \nabla^2 f({x}_1,{x}_2) \subseteq [\nabla^2f({x}_1,{x}_2)] = \begin{pmatrix}
    [0, 118] & [-69, 860] \\
    [-69, 860] & [0, 2152]
    \end{pmatrix}
    \label{eq:exampleIntHess}
\end{equation}
for $0 \leq {x}_1, {x}_2 \leq 2$. We also compute the interval Hessian by global minimization and maximization of each Hessian element using Gurobi and obtain the following: 
\begin{equation}
    [\nabla^2f({x}_1,{x}_2)]^* = \begin{pmatrix}
    [0, 118] & [-5, 860] \\
    [-5, 860] & [0, 2152]   
    \end{pmatrix}.
    \label{eq:exactIntHess}
\end{equation}
For this function, Gurobi-based optimization took approximately 47.1 milliseconds. In contrast, the interval arithmetic approach computed the bounds in roughly 0.0497 milliseconds, offering significantly faster but conservative bounds. Once the interval Hessian is computed, we estimate the minimum eigenvalue of the interval Hessian by using the following equation \cite{androulakis1995alphabb}: 
\begin{equation}
    \alpha \geq \max\left\{{0, -\frac{1}{2}\lambda_{min}([\nabla^2 f])}\right\},
    \label{eq:alpha1}
\end{equation}
where $\lambda_{min}$ is the minimum eigenvalue of the interval Hessian. Next, we discuss various methods to compute the minimum eigenvalue of an interval matrix. 
\subsection{\texorpdfstring{$\mathbf{\alpha}$ Calculation methods}{alpha Calculation methods}}
\label{sec:Methods}
\subsubsection{Approximate Gerschgorin theorem (GGN) for interval matrices}
This technique is based on applying the Gerschgorin theorem \cite{adjiman1996rigorous} on an interval matrix. For a real interval matrix $[ \mathbf{A} ] = ([\ubar{a}_{ij},\bar{a}_{ij}])$, the lower bound on the minimum eigenvalues is as follows:
\begin{equation}
    \lambda_{min}([\mathbf{A}]) \geq \min_i \bigg[\ubar{a}_{ii} - \sum_{j \neq i} \max ( |\ubar{a}_{ij}|,|\bar{a}_{ij}| ) \bigg].
    \label{eq:approxGGN}
\end{equation}
The computational complexity of this method is $\mathcal{O}(n^2)$. For example,  we use the nonlinear two variables function given in \eqref{eqn:example} to illustrate the use of \eqref{eq:approxGGN}. Using the interval Hessian in \eqref{eq:exampleIntHess} and \eqref{eq:exactIntHess}, we estimate $\lambda_{min}$ to be $-860$ for the two interval Hessian matrices.
\subsubsection{E-Matrix (EM) Method}
This method is an extension of the theorems developed by Deif \cite{deif1991interval} and Rohn \cite{rohn1998bounds} for the calculation of the lower bound on minimum eigenvalue $\lambda_{min}([\mathbf{A}])$ for a real interval matrix $[\mathbf{A}]$. In this method
\begin{equation}
    \lambda_{min}(\mathbf{[A]}) \geq \lambda_{min}(\widetilde{\mathbf{A}}_M + \mathbf{E}) - \rho(\widetilde{\Delta \mathbf{A}} + |\mathbf{E}|),
    \label{eq:EMatrix}
\end{equation}
where  $\Delta \mathbf{A} = (\Delta a_{ij})$ denotes a radius matrix with $\Delta a_{ij} = \dfrac{\overline{a}_{ij} - \underline{a}_{ij}}{2}$, $\widetilde{\Delta \mathbf{A}} = (\widetilde{\Delta a}_{ij})$ represents the modified radius matrix with
    \[
    \widetilde{\Delta a}_{ij} = 
    \begin{dcases}
    0 & \text{if } i = j, \\
    \Delta a_{ij} & \text{otherwise},
    \end{dcases}
    \]
 $\mathbf{A}_M = (a_{M,ij})$ denotes a mid-point matrix with $a_{M,ij} = \dfrac{\overline{a}_{ij} + \underline{a}_{ij}}{2}$,  $\tilde{\mathbf{A}}_M = (\widetilde{a}_{M,ij})$ represents the modified mid-point matrix with
    \[
    \widetilde{a}_{M,ij} = 
    \begin{dcases}
    \underline{a}_{ij} & \text{if } i = j, \\
    a_{M,ij} & \text{otherwise},
    \end{dcases}
    \]
and $\mathbf{E}$ is an arbitrary matrix that is taken to be equal to $\Delta \mathbf{A}$. The computational complexity of this method is $\mathcal{O}(n^3)$. We use the nonlinear two variables function given in \eqref{eqn:example} to illustrate the use of \eqref{eq:EMatrix}. Using the interval Hessian in \eqref{eq:exampleIntHess} and \eqref{eq:exactIntHess}, we estimate $\lambda_{min}$ to be $-1332.92$ and $-1331.88$, respectively.

\subsubsection{Mori-Kokame (MK) Method}
This method is an extension of the theorems developed by Mori and Kokame (MK) \cite{mori1994eigenvalue}, which suggest using lower and upper bound matrices to calculate the lower bound on the minimum eigenvalue for the interval matrix.
\begin{equation}
    \lambda_{min}([\mathbf{A}]) \geq \lambda_{min}(\underline{\mathbf{A}}) - \rho(\overline{\mathbf{A}} - \underline{\mathbf{A}}),
    \label{eq:MoriKokame}
\end{equation}
where $\overline{\mathbf{A}} = (\overline{a}_{ij})$ and $\underline{\mathbf{A}} = (\underline{a}_{ij})$.
The computational complexity of this method is $\mathcal{O}(n^3)$. We use the nonlinear two variables function given in \eqref{eqn:example} to illustrate the use of \eqref{eq:MoriKokame}. Using the interval Hessian in  \eqref{eq:exampleIntHess} and \eqref{eq:exactIntHess}, we estimate $\lambda_{min}$ to be $-2581.44$ and $-2475.11$, respectively. For this example, we observed that the GGN method, which is the least computationally expensive method, provides the tightest $\lambda_{min}$. However, this may not always be the case.  
\subsection{Algorithm Structure}
We incorporate the proposed search direction in a line-search framework and develop INTHOP - INTerval Hessian-based OPtimization. In this section, we present two variants of INTHOP based on fixed and adaptive interval size selection, represented as INTHOP:F and INTHOP:A, respectively.
\subsubsection{INTHOP:F}
The line search framework finds the next iterate $\mathbf{x}_{k+1}$ based on the search direction $\mathbf{p}_k$ and the step length $\theta_k$
\begin{equation}
    \mathbf{x}_{k+1} = \mathbf{x}_k+\theta_k\mathbf{p}_k.
\end{equation}

The INTHOP algorithm with fixed interval size, denoted as INTHOP:F, is given in Algorithm~\ref{algo:GlobalStructure}. We select two iteration counters, $t$ and $k$. The former functions as the outside iteration counter, while the latter serves as the inner iteration counter.  The iteration counter $t$ is updated whenever Hessian is computed. In contrast, $k$ is updated whenever a new iterate is generated. At the beginning of iteration $k$, we check whether the gradient norm is less than the pre-defined tolerance or if the number of iterations has exceeded. If either of the conditions is satisfied, we terminate and return the solution. Otherwise, we check whether the current iterate $\mathbf{x}_k$ lies in the interval $[\mathbf{x}_t^L, \mathbf{x}_t^U]$. To avoid confusion, we note that $\mathbf{x}_t$ is equal to that  iterative point $\mathbf{x}_{k}$  for which Hessian is updated whenever $\mathbf{x}_{k}$ does not belong to the interval $[\mathbf{x}_{t-1}^L, \mathbf{x}_{t-1}^U]$.

The Hessian modification becomes essential whenever the current iterate \(\mathbf{x}_k\) exits the interval around which the previous Hessian was computed because the difference between $\nabla^2 \mathcal{H}_t$ and $\nabla^2 {f}_k$ may not be bounded. At such points, we establish a new interval of length \(\delta\) centered at \(\mathbf{x}_t\), and compute the interval Hessian within this region. This localized approximation captures the curvature of the objective function over the updated interval and informs the subsequent search directions. One approach to construct intervals is to expand the interval around the current iterate symmetrically, as given in Algorithm \ref{algo:SymInt}. 

Once the interval is constructed, the interval Hessian and the corresponding $\alpha$ are calculated using one of the methods outlined in Section \ref{sec:Methods}. Thereafter, we compute $\nabla^2\mathcal{H}_t$, its inverse, and the search direction. The search direction $\mathbf{p}_k$ fulfills ($\mathbf{g}_k^{\top}\mathbf{p}_{k}<0$), ensuring a reduction in the function value. Consequently, $f_{k+1} < f_k$ after each iteration. The step length is found by backtracking line search, given in Algorithm \ref{algo:LineSearch}, such that the Armijo condition is satisfied. The new iterate is then obtained by taking a step length $\theta_k$ in the direction $\mathbf{p}_k$ from the current iterate $\mathbf{x}_k$. The procedure is repeated until a convergence criterion is satisfied.  

\begin{figure}[H]
\centering
\includegraphics[width=0.5\textwidth]{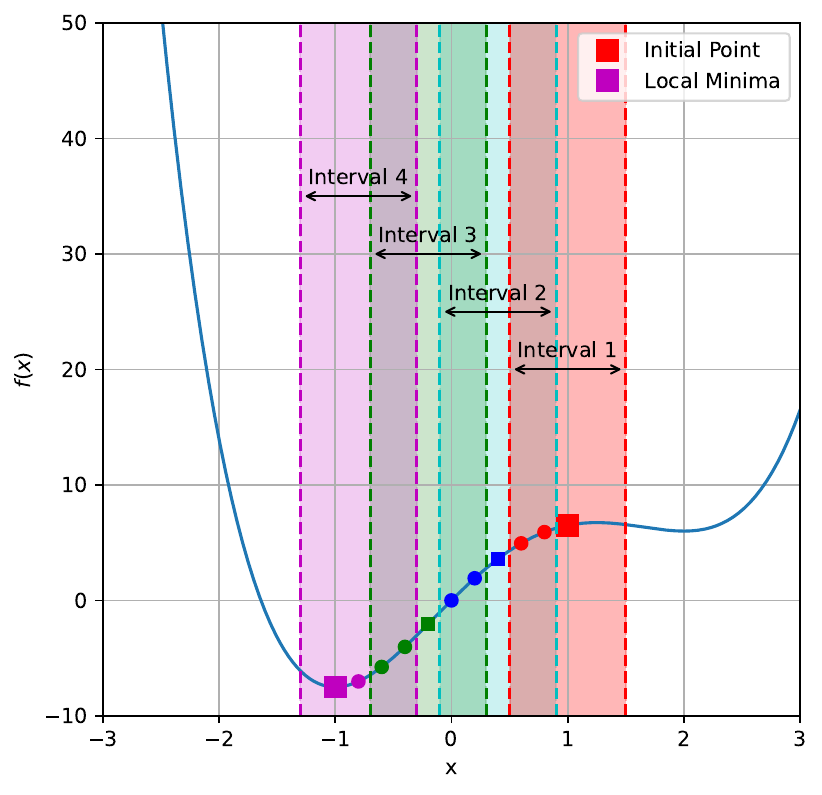}
\caption{Graphical illustration of INTHOP with fixed interval size (INTHOP:F) and Algorithm \ref{algo:SymInt} used to generate intervals.}
\label{fig:graphical_representation of algorithm}
\end{figure}

Figure \ref{fig:graphical_representation of algorithm} demonstrates the iterative mechanics of our algorithm using a one-dimensional function with the bounds updating scheme given in Algorithm \ref{algo:SymInt}. The figure shows that a new interval of length $\delta$, centered around the current iterate, represented by the square markers, is created. The new interval is created at the starting point or when the current iterate $\mathbf{x}_k \notin [\mathbf{x}_t^L, \mathbf{x}_t^U]$. At these points, the expensive computations occur. First, the interval Hessian is computed. Second, the lower bound on the minimum eigenvalue for the interval Hessian is estimated, and the Hessian matrix is modified to ensure positive definiteness. Third, the modified Hessian matrix is inverted to compute the search direction. Circle markers represent points where $\mathbf{x}_k \in [\mathbf{x}_t^L, \mathbf{x}_t^U]$. In these cases, only the gradient is evaluated, and the search direction is computed by multiplying the stored inverse Hessian from the $t^{\text{th}}$ iteration by the gradient vector.

\begin{algorithm}[H]
\caption{INTHOP:F}
\textit{Init.} Choose $\mathbf{x}_0 \in \mathbb{R}^n,\ \epsilon_g \in (0, 1),\ t = 0, c_1>0;$\par
\begin{algorithmic}
\For {$k = 0,1,2,\dots$}
    \If {$\|\mathbf{g}_k\| < \epsilon_g$}
        \State \textbf{Terminate - Solution found};
        \vspace{0.3em}
    \EndIf \textbf{end if}
    \If {$k>iter_{max}$}
        \State \textbf{Terminate}
        \vspace{0.3em}
    \EndIf \textbf{end if}
    \If {$k=0$ \text{or} $\mathbf{x}_k \notin [\mathbf{x}_t^L,\mathbf{x}_t^U]$}  
        \State $t \xleftarrow{} t + 1$
        \State Go to \textbf{Algorithm~\ref{algo:SymInt}} to compute $\mathbf{x}^L_t$ and $\mathbf{x}^U_t$
        \State Choose  a constant $\tilde g >0$ such that $\tilde g \geq \lVert \textbf{g}_{t}\rVert$
        \State Go to \textbf{Algorithm~\ref{algo:HessianMod}} to compute $\nabla^2 \mathcal{H}_t$
         
    \vspace{0.5em}
    \EndIf \textbf{end if}
    
   \State Find the search direction by solving
    \begin{equation}\label{p_k}
        \mathbf{p}_k = - \nabla^2 \mathcal{H}_t^{-1}\mathbf{g}_k
    \end{equation}
   \State Compute step length $\theta_k$ using Algorithm \ref{algo:LineSearch} and set $\mathbf{x}_{k+1} = \mathbf{x}_k + \theta_k \mathbf{p}_k$
\EndFor \textbf{end for}
\end{algorithmic}
\label{algo:GlobalStructure}
\end{algorithm}


\begin{algorithm}[H]
\caption{Symmetric interval generation scheme}
\textit{Init. }$\delta \in \mathbb{R}^+$
\begin{algorithmic}
\For {$i = 0,1,2,\dots, n$}
    \State $\mathbf{x}^L_{t,i} = \mathbf{x}_{k,i} - \delta / 2$
    \vspace{0.2em}
    \State $\mathbf{x}^U_{t,i} = \mathbf{x}_{k,i} + \delta / 2$
    \vspace{0.3em}
\EndFor \textbf{end for}
\end{algorithmic}
\label{algo:SymInt}
\end{algorithm}

\begin{algorithm}[H]
\caption{Hessian approximation algorithm}
\begin{algorithmic}
\State Compute the Interval Hessian $[\nabla^2 f(\mathbf{x})]$ for $\mathbf{x} \in [\mathbf{x}_t^L,\mathbf{x}_t^U] $ using interval arithmetic
\State Compute $\lambda_{t, min}$ using Equation \ref{eq:approxGGN}, \ref{eq:EMatrix} or \ref{eq:MoriKokame}
\begin{equation}
    \alpha = \max \left(0, -\frac{1}{2} \lambda_{t, min} \right)
\end{equation}
\State Compute the positive semidefinite Hessian approximation,
\begin{equation}
    \nabla^2 \mathcal{H}_t = \nabla^2f_t + 2\alpha I+c_{1}\tilde{g} I
\end{equation}
\end{algorithmic}
\label{algo:HessianMod}
\end{algorithm}

\begin{algorithm}[H]
\caption{Backtracking line search}
Choose $\theta_0 > 0,\space \eta\in (0,1),\space \nu \in (0, 1)$ \\
$\theta \gets \theta_0$
\begin{algorithmic}
\While{$f(\mathbf{x}_{k} + \theta\mathbf{p}_k) > f (\mathbf{x}_k) + \eta\theta\nabla f(\mathbf{x}_k)^\top\mathbf{p}_k$}
    \State $\theta \gets \nu \theta$
\EndWhile
\State \textbf{end while} \newline
$\theta_k \gets \theta$ \newline
$\mathbf{x}_{k+1} \gets \mathbf{x}_{k} + \theta_k\mathbf{p}_k$
\end{algorithmic}
\label{algo:LineSearch}
\end{algorithm}

\subsubsection{INTHOP:A} \label{sec:INTHOP:A}
INTHOP:F (Algorithm \ref{algo:GlobalStructure}) assumed a constant value of the interval size, making it a hyperparameter for our algorithm. However, the same value of interval size might not be appropriate for a different function or different regions of a function. Accordingly, we develop two variants of INTHOP based on adaptive interval size.

\newtheorem{lemma}{Lemma}


A high $\alpha$ value ensures that $\nabla^2\mathcal{H}_t$ is positive definite and, consequently, the search direction is descent. However, it leads to conservative steps and potentially slow convergence. Therefore, an appropriate balance is required to ensure positive definiteness and step lengths. We develop two INTHOP variants based on adaptively updating the interval size. 

\subsubsection{INTHOP:A1}

This method is based on the observation that the magnitude of the search direction can serve as an indicator of the quality of the Hessian approximation. Specifically, a small search direction norm may imply that the algorithm is approaching a stationary point (small gradient norm) or that $\alpha$ is too large due to a large interval size. Conversely, a large search direction indicates that the interval is sufficiently small and the Hessian approximation is reliable.

To quantify this, we define the following ratio:

\begin{equation}\label{eq:A1}
    \tau_t = \frac{r}{\sqrt{n}} \times \frac{\| \mathbf{p}_{k-1} \|_1}{\sqrt{\| \mathbf{p}_{k-1} \|_2^2 + \beta}},
\end{equation}
\noindent
where \( \beta > 0 \) is a small regularization constant to avoid division by zero and dampening the ratio $\tau_t$. The parameter \( r \) is a scaling factor that adjusts the interval size. This ratio $\tau_t$ is inspired by the adaptive learning rate strategies used by machine learning optimizers such as Adam \cite{kingma2014adam}. These methods combine first- and second-moment estimates of the gradients to modulate the step size adaptively to improve convergence rates. Adam, for instance, adjusts the step size by dividing the running average of the gradient (first moment) by the square root of the running average of the squared gradient (second moment) with a small constant added to prevent division by zero. This dampening term ($\beta$) also smooths the update to prevent instability from sudden gradient spikes.

For any nonzero vector \( \mathbf{p}_k \in \mathbb{R}^n \), the ratio is bounded as:

\begin{equation}
    0 \leq \tau_t \leq r \frac{\| \mathbf{p}_{k-1} \|_1}{\sqrt{\| \mathbf{p}_{k-1} \|_2^2 + \beta}}.
\end{equation}
\noindent
In the regime where \( \| \mathbf{p}_k \|_2^2 \gg \beta \), the upper bound asymptotically approaches r. When the search direction is small in magnitude, the ratio \( \tau_t \) is small, indicating potential over-regularization and suggesting a reduction in the interval size. Conversely, a higher \( \tau_t \) suggests that the current interval size is effective, and it can be increased to allow larger steps. The INTHOP variant with adaptive interval size based on Eq. \ref{eq:A1} is denoted as INTHOP:A1 and summarized in Algorithm~\ref{algo:AdaptGlobalStructure}.
\begin{algorithm}[H]
\caption{INTHOP:A1}
\textit{Init.} Choose $\mathbf{x}_0 \in \mathbb{R}^n,\ \epsilon_g \in (0, 1),\ \delta = \delta_0,\ r > 0,\ \beta > 0,\ t = 0, c_1>0;$\par
\begin{algorithmic}
\For {$k = 0,1,2,\dots$}
    \If {$\|\mathbf{g}_k\| < \epsilon_g$}
        \State \textbf{Terminate - Solution found};
        \vspace{0.3em}
    \EndIf \textbf{end if}
    \If {$k>iter_{max}$}
        \State \textbf{Terminate}
        \vspace{0.3em}
    \EndIf \textbf{end if}
    \If {$k=0$ \text{or} $\mathbf{x}_k \notin [\mathbf{x}_t^L,\mathbf{x}_t^U]$}
        \If{$k=0$}
            \State $\delta_t = \delta_0$
        \Else
            \State  Calculate
                \begin{equation} 
                    \tau_t = \frac{r}{\sqrt{n}} \times \frac{\| \mathbf{p}_{k-1} \|_1}{\sqrt{\| \mathbf{p}_{k-1} \|_2^2 + \beta}}
                \end{equation}
                Update
                \begin{equation}
                    \delta_t = \delta_{t-1}\times \tau_t
                \end{equation}
        \EndIf \textbf{end if}
        \State $t \xleftarrow{} t + 1$
        \State Go to \textbf{Algorithm~\ref{algo:SymInt}} to compute $\mathbf{x}^L_t$ and $\mathbf{x}^U_t$
        \State Choose  a constant $\tilde g >0$ such that $\tilde g \geq \lVert \textbf{g}_{t}\rVert$
        \State Go to \textbf{Algorithm~\ref{algo:HessianMod}} to compute $\nabla^2 \mathcal{H}_t$
        \State Compute $\nabla^2 \mathcal{H}_t^{-1}$
    \vspace{0.5em}
    \EndIf \textbf{end if}
    
    \State Find the search direction by solving
    \begin{equation}
        \mathbf{p}_k = - \nabla^2 \mathcal{H}_t^{-1}\mathbf{g}_k
    \end{equation}
    \State Compute step length $\theta_k$ using Algorithm \ref{algo:LineSearch} and set $\mathbf{x}_{k+1} = \mathbf{x}_k + \theta_k \mathbf{p}_k$
    \vspace{0.5em}
\EndFor \textbf{end for}
\end{algorithmic}
\label{algo:AdaptGlobalStructure}
\end{algorithm}
\subsubsection{INTHOP:A2}
We created a modified second-order Taylor model (Equation~\ref{eq:Modified-Model}) to obtain the search direction within a given interval. INTHOP algorithm can also be interpreted through the lens of the trust-region method as follows. From the \( t^{\text{th}} \) to the \( (t+1)^{\text{th}} \) iterate, we move within a fixed-size interval constructed around the iterate \( \mathbf{x}_t \). Once an iterate reaches any boundary of this interval, a new interval is formed. This behavior is similar to classical trust-region methods, where for a given local model, we define a region within which we trust the model to approximate the objective function. The size of this region is important. If it is too small, the model may be accurate, but progress is slow; if it is too large, the model (Eq. \ref{eq:Modified-Model}) may not accurately represent the true function. In trust-region methods, a
ratio is computed between the actual reduction and the predicted reduction (i.e., the decrease predicted by the local model), and this ratio is used to adjust the trust-region size accordingly.

In our case, we update the interval size based on the performance of the algorithm in the previous interval. We define the ratio $\xi_{t+1}$ here as follows:
\begin{equation}\label{eq:A2}
    \xi_{t+1} = \frac{f(\mathbf{x}_{t}) - f(\mathbf{x}_{t+1})}{-(\mathbf{g}_t^\top \mathbf{s}_t + \mathbf{s}_t^\top \nabla^2 \mathcal{L}_{t}\mathbf{s}_t)},
\end{equation}
where $\mathbf{s}_t = \mathbf{x}_{t+1} - \mathbf{x}_t$. Here, the numerator is the decrease in the actual function value in the last interval, and the denominator shows the decrease obtained by the model. The following three cases may arise: (1) if $\tau_t < 0$, then the objective function increases while the model predicted it to decrease, which means the previous value of the interval size gave an inaccurate representation of the objective function. (2) If $\tau_t$ is close to zero and positive, then the model has significantly overestimated the decrease in the objective function. (3) If \( \tau_t \) is close to 1, the model's predicted decrease closely matches the actual decrease in the objective function. In the first two cases, the interval size is reduced to improve the accuracy of the Hessian approximation. For the third case, the current interval size is retained or increased depending on the value of $\tau_t$. The INTHOP variant based on \eqref{eq:A2} to update the interval size is given in Algorithm \ref{algo:TrustAdaptGlobalStructure}. 
\begin{algorithm}[H]
\caption{INTHOP:A2}
\textit{Init.} Choose $\mathbf{x}_0 \in \mathbb{R}^n,\ \epsilon_g \in (0, 1),\ \delta = \delta_0,\ iter_{max} > 0,\ \delta_{max} > \delta_{min} > 0,\ t = 0, c_1>0;$\par
\begin{algorithmic}
\For {$k = 0,1,2,\dots$}
    \If {$\|\mathbf{g}_k\| < \epsilon_g$}
        \State \textbf{Terminate - Solution found};
        \vspace{0.3em}
    \EndIf \textbf{end if}
    \If {$k>iter_{max}$}
        \State \textbf{Terminate}
        \vspace{0.3em}
    \EndIf \textbf{end if}
    \If {$k=0$ \text{or} $\mathbf{x}_k \notin [\mathbf{x}_t^L,\mathbf{x}_t^U]$}
        \If{$k=0$}
            \State $\delta_t = \delta_0$
        \Else        
            \State $t \xleftarrow{} t + 1$
            \State Compute $\xi_t$ with
                \begin{equation}
                    \xi_t = \frac{f(\mathbf{x}_{t-1}) - f(\mathbf{x}_{t})}{-(\mathbf{g}_{t-1}^\top \mathbf{s}_{t-1} + \mathbf{s}_{t-1}^\top \nabla^2 \mathcal{L}_{t-1}\mathbf{s}_{t-1})}
                \end{equation}
                \If{$\xi_t < \frac{1}{4}$}
                    \State $\delta_{t} \gets \max\left(\frac{\delta_{t-1}}{2}, \delta_{\min}\right)$
                \ElsIf{$\xi_t > \frac{3}{4}$}
                    \State $\delta_{t} \gets \min(4 \cdot \delta_{t-1}, \delta_{\max})$
                \ElsIf{$\frac{1}{4} \leq \xi_t \leq \frac{3}{4}$}
                        \State $\delta_{t} \gets \delta_{t-1}$
                    \vspace{0.5em}
                \EndIf \textbf{end if}
                \vspace{0.5em}
        \EndIf \textbf{end if}
        \State Go to \textbf{Algorithm~\ref{algo:SymInt}} to compute $\mathbf{x}^L_t$ and $\mathbf{x}^U_t$
        \State Choose  a constant $\tilde g >0$ such that $\tilde g \geq \lVert \textbf{g}_{t}\rVert$
        \State Go to \textbf{Algorithm~\ref{algo:HessianMod}} to compute $\nabla^2 \mathcal{H}_t$
        \State Compute $\nabla^2 \mathcal{H}_t^{-1}$
    \vspace{0.5em}
    \EndIf \textbf{end if}
    
    \State Find the search direction by solving
    \begin{equation}
        \mathbf{p}_k = - \nabla^2 \mathcal{H}_t^{-1}\mathbf{g}_k
    \end{equation}
    \State Compute step length $\theta_k$ using Algorithm \ref{algo:LineSearch} and set $\mathbf{x}_{k+1} = \mathbf{x}_k + \theta_k \mathbf{p}_k$
    \vspace{0.5em}
\EndFor \textbf{end for}
\end{algorithmic}
\label{algo:TrustAdaptGlobalStructure}
\end{algorithm}

{\section{Global Convergence}
In this section, we prove the convergence of the proposed method.  First, we prove important inequalities related to the bounds on the direction $\mathbf{p}_{k}$ and the step length calculated by Armijo line search. Then, the global convergence of the proposed method is established in the sense that there exists an accumulation point of the sequence generated by the algorithm that is a stationary point.

Next, we mention an assumption that is considered to prove the following results in this section.
 

\begin{assumption}
    For every sequence $\{\mathbf{x}_{k}\}$ generated by the algorithm, there exists a compact set $\mathcal{S}\subseteq\mathbb{R}^{n}$ that contains sequence $\{\mathbf{x}_{k}\}$, i.e., $\{\mathbf{x}_{k}\}\subseteq\mathcal{S}$.
\end{assumption}
    From now on, we assume that $\lVert \mathbf{g}_k\rVert \neq 0$ for all $k\geq0$. It means that there exists an $\epsilon>0$ such that 
\begin{equation}\label{normgklessepsilon}
\lVert \mathbf{g}_k\rVert \geq\epsilon \text{ for all } k\geq0.
\end{equation}
Also, note that $\{\mathbf{x}_{k}\}$ belongs to the compact set $\mathcal{S}$ and $f$ is continuously differentiable. Therefore, there exists a  constant $\mathcal{G}$ such that
\begin{equation}\label{normgklessG}
\lVert \mathbf{g}_k\rVert \leq \mathcal{G} \text{ for all } k\geq0.
\end{equation}


\begin{lemma}\label{normpklessc1gk}
    Suppose that $\{\mathbf{p}_{k}\}$ is the sequence of directions generated by the algorithm. Then, 
    \[
    \lVert \mathbf{p}_{k}\rVert \leq  \frac{\mathcal{G}}{ c_1\tilde g}.
    \]
\end{lemma}
\begin{proof}
Note that $\nabla^{2}\mathcal{H}^{-1}_{t}=(\nabla^{2}\mathcal{L}_{t}+c_1\tilde g I)^{-1}$, where $\nabla^{2}\mathcal{L}_{t}=\nabla^{2}\mathcal{L}(\mathbf{x}_{t})=\nabla^{2} f(\mathbf{x}_{t})+2\alpha I$ and $\mathbf{x}_{t}\in
    [\mathbf{x}^{L}_{t}, \mathbf{x}^{U}_{t}]$. Therefore, 
    \begin{align}
    \lVert (\nabla^{2} f(\mathbf{x}_{t})+2\alpha I+c_1\tilde g I)^{-1}\rVert
        &=\lambda_{\max} \left((\nabla^{2} f(\mathbf{x}_{t})+2\alpha I+c_1\tilde g I)^{-1}\right)\nonumber\\
        &=\frac{1}{\lambda_{\min} \left(\nabla^{2} f(\mathbf{x}_{t})+2\alpha I+c_1\tilde g I\right)}\nonumber\\
        & \leq \frac{1}{\lambda_{\min} \left(\nabla^{2} f(\mathbf{x}_{t})+2\alpha I\right)+c_1\tilde g } \label{lemma:regularization_effect}\\
        &\leq \frac{1}{c_1\tilde g } \text{ as } \lambda_{\min} \left(\nabla^{2} f(\mathbf{x}_{t})+2\alpha I\right)\geq 0\label{normHkless}.
    \end{align}
    For each $k\geq0$, we have $\mathbf{p}_{k}=-\nabla^{2}\mathcal{H}^{-1}_{t}\mathbf{g}_k=-(\nabla^{2}\mathcal{L}_{t}+c_1\tilde g I)^{-1}\mathbf{g}_k$  which implies that 
    \begin{align}\label{pklessinvres}
        \lVert \mathbf{p}_{k}\rVert & =\lVert (\nabla^{2} f(\mathbf{x}_{t})+2\alpha I+c_1\tilde g I)^{-1} \mathbf{g}_k\rVert\nonumber\\
        & \leq   \lVert (\nabla^{2} f(\mathbf{x}_{t})+2\alpha I+c_1\tilde g I)^{-1}\rVert \lVert \mathbf{g}_k\rVert \\ 
        & \leq \frac{\lVert g_k\rVert}{c_1\tilde g } \text{ using } \eqref{normHkless}\nonumber\\
        & \leq \frac{\mathcal{G}}{ c_1\tilde g}\text{ using } \eqref{normgklessG}\nonumber,
    \end{align}
    and the proof is complete.
\end{proof}

\begin{remark}
    From Lemma \ref{normpklessc1gk}, we have that $\mathbf{x}_{k}+\theta \mathbf{p}_{k}\in \mathcal{S}+\mathbb{B}(0,\mathcal{D})$, where $\theta\in(0,1]$ and $\mathbb{B}$ is a closed ball centered at $0$ and with radius $\mathcal{D}=\frac{\mathcal{G}}{ c_1\tilde g}$. Note that $\mathcal{S}+\mathbb{B}(0,\mathcal{D})$ is a compact set and $f$ is  a twice continuously differentiable function. Therefore, there exists a constant $\mathcal{A}$ such that 
    \begin{equation}\label{steplengthpf3}
        \lVert \nabla^{2} f(x)\rVert \leq \mathcal{A} \text{ for all } x\in \mathcal{S}+\mathbb{B}(0,\mathcal{D}).
    \end{equation}
\end{remark}

In the next result, it is proved that the step length, calculated using Armijo line search condition, has a positive lower bound.

\begin{lemma}\label{steplengthlemma}
    Suppose that $\theta_{k}$ is the step length generated by the algorithm using Armijo line search condition for each $k\geq0$. Then,  
    \begin{equation*}
    \theta_{k}\geq \min\left(1,\frac{2(1-\eta)c_1\tilde{g}}{\mathcal{A}}{t}\right),  \text{ where } 0<t<1 \text{ and } \mathcal{A} \text{ is given in }\eqref{steplengthpf3}.      
    \end{equation*}
\end{lemma}
\begin{proof}
    Note that $f$ is twice differentiable function. Therefore, from Taylor's theorem, there exists a scalar $b\in (0,1)$ such that 
    \begin{align}\label{steplengthpf1}
    &f(\mathbf{x}_{k}+\theta_{k}\mathbf{p}_k)=-f(\mathbf{x}_{k})+\theta_{k}\mathbf{g}^{\top}_{k}\mathbf{p}_{k}+\frac{1}{2} \theta^{2}_{k}\mathbf{p}^{\top}_{k} \nabla^{2} f(\mathbf{x}_{k}+b\theta_{k}\mathbf{p}_{k})\mathbf{p}_{k}\nonumber\\
    \implies& f(\mathbf{x}_{k})-f(\mathbf{x}_{k}+\theta_{k}\mathbf{p}_{k})+\eta\theta_{k}\mathbf{g}^{\top}_{k}\mathbf{p}_{k}=-(1-\eta)\theta_{k}\mathbf{g}^{\top}_{k}\mathbf{p}_{k}-\frac{1}{2} \theta^{2}_{k}\mathbf{p}^{\top}_{k}\nabla^{2}f(\mathbf{x}_{k}+b\theta_{k}\mathbf{p}_{k})\mathbf{p}_{k}.
    \end{align}
    From \eqref{p_k}, we have $\mathbf{p}_{k}=-\nabla^{2}\mathcal{H}^{-1}_{t}\mathbf{g}_k$, where  $\nabla^{2}\mathcal{H}^{-1}_{t}=(\nabla^{2} f(\mathbf{x}_{t})+2\alpha I+c_1\tilde g I)^{-1}$ and $\mathbf{x}_{t}\in
    [\mathbf{x}^{L}_{t}, \mathbf{x}^{U}_{t}]$. Therefore,
    \begin{equation}\label{steplengthpf2}
    \mathbf{g}_k=-(\nabla^{2} f(\mathbf{x}_{t})+2\alpha I+c_1\tilde g I)\mathbf{p}_{k}.
    \end{equation}
    In view of \eqref{steplengthpf1} and \eqref{steplengthpf2}, we get 
    \begin{align}\label{steplengthpf4}
        &~~~~f(\mathbf{x}_{k})-f(\mathbf{x}_{k}+\theta_{k}\mathbf{p}_{k})+\eta\theta_{k}\mathbf{g}^{\top}_{k}d_{k}\nonumber\\
        &= (1-\eta) \theta_{k} \mathbf{p}^{\top}_{k} (\nabla^{2} f(\mathbf{x}_{t})+2\alpha I+c_1\tilde g I)\mathbf{p}_{k}-\frac{1}{2} \theta^{2}_{k} \mathbf{p}^{\top}_{k}\nabla^{2} f(\mathbf{x}_{k}+b\theta_{k}\mathbf{p}_{k})\mathbf{p}_k\nonumber\\
        &= (1-\eta)\theta_{k}\mathbf{p}^{\top}_{k}(\nabla^{2} f(\mathbf{x}_{t})+2\alpha I)\mathbf{p}_{k}+(1-\eta)\theta_{k}\mathbf{p}^{\top}_{k}\left(c_1\tilde gI-\frac{1}{2(1-\eta)}\theta_{k}\nabla^{2} f(\mathbf{x}_{k}+b\theta_{k}\mathbf{p}_{k})\right)\mathbf{p}_{k}\nonumber\\
        &\geq (1-\eta) \theta_{k} \mathbf{p}^{\top}_{k} \left(c_1\tilde gI-\frac{1}{2(1-\eta)}\theta_{k}\nabla^{2} f(\mathbf{x}_{k}+b\theta_{k}\mathbf{p}_{k})\right)\mathbf{p}_{k}\nonumber\\
        &~~~~~~~~~~~~~~~~~~~~~~~~~~~~~~~~~~~~~~~~~~~~~~~\text{ as }\nabla^{2} f(\mathbf{x}_{t})+2\alpha I \text{ is positive semidefinite  for }\mathbf{x}_{t}\in[\mathbf{x}^{L}_{t}, \mathbf{x}^{U}_{t}]\nonumber\\
        &=(1-\eta)\theta_{k}\left(c_1\tilde g\lVert \mathbf{p}_{k}\rVert^{2}I-\frac{1}{2(1-\eta)}\theta_{k}\mathbf{p}^{\top}_{k}\nabla^{2}f(\mathbf{x}_{k}+b\theta_{k}\mathbf{p}_{k})\mathbf{p}_{k}\right)\nonumber\\
        &\geq (1-\eta)\theta_{k}\left(c_1\tilde g-\frac{1}{2(1-\eta)}\theta_{k}\rVert\nabla^{2}f(\mathbf{x}_{k}+b\theta_{k}\mathbf{p}_{k})\rVert\right)\lVert \mathbf{p}_{k}\rVert^{2}.
    \end{align}
    
    Now, from  \eqref{steplengthpf4}, and \eqref{steplengthpf3}, we get 
    \begin{align}\label{steplengthpf5}
        f(\mathbf{x}_{k})-f(\mathbf{x}_{k}+\theta_{k}\mathbf{p}_{k})+\eta\theta_{k}\mathbf{g}^{\top}_{k}d_{k}\geq (1-\eta) \theta_{k} \left(c_1\tilde{g}-\frac{1}{2(1-\eta)}\theta_{k}\mathcal{A}\right)\lVert \mathbf{p}_{k}\rVert^{2}.
    \end{align}
Next, we consider two cases.\\
First, if $\frac{2(1-\eta)c_1\tilde{g}}{\mathcal{A}}\geq1$, then from \eqref{steplengthpf5}, 
\begin{equation*}
f(\mathbf{x}_{k})-f(\mathbf{x}_{k}+\mathbf{p}_k)\geq-\eta \mathbf{g}^{\top}_{k}\mathbf{p}_{k}
\end{equation*}
which implies that $\theta_{k}=1$ satisfies Armijo line search condition.\\
Second, if $\frac{2(1-\eta)c_1\tilde{g}}{\mathcal{A}}<1$, then for $\theta_{k}\leq \frac{2(1-\eta)c_1\tilde{g}}{\mathcal{A}}$ in \eqref{steplengthpf5}, we have 
\begin{equation}
    f(\mathbf{x}_{k})-f(\mathbf{x}_{k}+\mathbf{p}_k)\geq-\eta \mathbf{g}^{\top}_{k}\mathbf{p}_{k}
\end{equation}
which implies that $\theta_{k}$ must be greater than or equal to $\frac{2(1-\eta)c_1\tilde{g}}{\mathcal{A}}{t_1}$ for $0<t_1<1$. Otherwise, there will be a $\tilde{\theta_{k}}=\frac{\theta_{k}}{t_1}>\theta_{k}$ that 
satisfies 
$ f(\mathbf{x}_{k})-f(\mathbf{x}_{k}+\tilde{\theta_{k}}\mathbf{p}_k)\geq-\eta\tilde{\theta_{k}} \mathbf{g}^{\top}_{k}\mathbf{p}_{k}$ which contradicts the definition of $\theta_{k}$. Hence, we have $\theta_{k}\geq\min\left(1,\frac{2(1-\eta)c_1\tilde{g}}{\mathcal{A}}{t_1}\right)$, where  $0<t_1<1$. 
\end{proof}
Next, it is proved that the difference between the value of the function on two iterative points generated by the algorithm has a lower bound.
\begin{lemma}\label{differencelemma}
    Suppose that the sequence $\{\mathbf{x}_{k}\}$ is generated by the algorithm. Then,
\[
f(\mathbf{x}_{k})-f(\mathbf{x}_{k+1})\geq \frac{\eta\epsilon^{2}\theta_{\min}}{(1+b)\mathcal{A}+c_1\tilde{g}},
\]
where $\theta_{\min}= \min\left(1,\frac{2(1-\eta)c_1\tilde{g}}{\mathcal{A}}{t}\right)$ with $0<t<1$.
\end{lemma}
\begin{proof}
    From Armijo line search condition, we have 
    \begin{align}\label{differencelemmapf1}
        f(\mathbf{x}_{k})-f(\mathbf{x}_{k+1})&\geq -\eta\theta_{k}\mathbf{g}^{\top}_{k}\mathbf{p}_{k}\nonumber\\
                        & = \eta\theta_{k}\mathbf{g}^{\top}_{k}(\nabla^{2} f(\mathbf{x}_{t})+2\alpha I+c_1\tilde gI)^{-1}\mathbf{g}_k\text{ using the definition of }\mathbf{p}_{k}\nonumber\\
                        &\geq \eta \theta_{k}  \lambda_{\min}\left((\nabla^{2} f(\mathbf{x}_{t})+2\alpha I+c_1\tilde gI)^{-1}\right)\lVert \mathbf{g}_k\rVert^{2}.
    \end{align}
    It can be easily note that
    \begin{align}\label{differencelemmapf2}
        \lambda_{\min}\left((\nabla^{2} f(\mathbf{x}_{t})+2\alpha I+c_1\tilde gI)^{-1}\right)&=\frac{1}{\lambda_{\max}\left(\nabla^{2} f(\mathbf{x}_{t})+2\alpha I+c_1\tilde gI\right)}\nonumber\\
        &=\frac{1}{\lambda_{\max}\left(\nabla^{2} f(\mathbf{x}_{t})\right)+2\alpha+c_1\tilde g}\nonumber\\
        &=\frac{1}{\lVert\nabla^{2} f(\mathbf{x}_{t})\rVert+2\alpha+c_1\tilde g}\nonumber\\
        &\geq \frac{1}{(\mathcal{A}+2\alpha+c_1\tilde{g})} \text{ using }\eqref{steplengthpf3}.
    \end{align}
    Thus, form \eqref{differencelemmapf1} and \eqref{differencelemmapf2}, we have 
    \begin{align*}
        f(\mathbf{x}_{k})-f(\mathbf{x}_{k+1})&\geq \eta\theta_{k}\frac{1}{(\mathcal{A}+2\alpha+c_1\tilde{g})}\lVert \mathbf{g}_k\rVert^{2}\\
        &\geq \eta\theta_{\min}\frac{1}{(\mathcal{A}+2\alpha+c_1\tilde{g})}\epsilon^{2} \text{ using Lemma }\ref{steplengthlemma} \text{ and } \eqref{normgklessepsilon},
    \end{align*}
    where $\theta_{\min}= \min\left(1,\frac{2(1-\eta)c_1\tilde{g}}{\mathcal{A}}{t}\right)$ with $0<t<1$, and hence the proof is complete.
\end{proof}

Finally, using the results proved in this section, we prove the global convergence of the algorithm.
\allowdisplaybreaks
\begin{theorem}
    Suppose that the sequence $\{\mathbf{x}_{k}\}$ is generated by the algorithm. Then, 
    \begin{equation*}
        \lim\limits_{k\to\infty}\lVert \mathbf{g}_k\rVert =0.
    \end{equation*}
    \begin{proof}
        We shall prove this theorem by the method of contradiction. Let the limit
        $\lim\sup\limits_{k\to\infty}\lVert \mathbf{g}_k\rVert\neq0$. Also, suppose that $\epsilon=\frac{\lim\sup\limits_{k\to\infty}\lVert \mathbf{g}_k\rVert}{2}$ and $\mathcal{J}_{k,\epsilon}=\{j\in\{0,1,2,\ldots\} : j\leq k, \lVert g_{j}\rVert \geq\epsilon\}$. Moreover, we assume that $\omega_{k}$ is the cardinality (number of elements in the set) of the set $\mathcal{J}_{k,\epsilon}$. Then, we have 
        \begin{equation}\label{cgspf1}
            \lim\limits_{k\to\infty}\omega_k=\infty.
        \end{equation}
        Note that
        \begin{align*}
            f(\mathbf{x}_{0})-f(\mathbf{x}_{k+1})&=\sum\limits_{j=0}^{k}(f(\mathbf{x}_{j})-f(\mathbf{x}_{j+1}))\\
            &\geq \sum\limits_{j\in\mathcal{J}_{k,\epsilon}}(f(\mathbf{x}_{j})-f(\mathbf{x}_{j+1}))\\
            &\geq \sum\limits_{j\in\mathcal{J}_{k,\epsilon}} \left(\frac{\eta\epsilon^{2}\theta_{\min}}{(1+b)\mathcal{A}+c_1\tilde{g}}\right) \text{ using Lemma }\ref{differencelemma}\\
            &= \left(\frac{\eta\epsilon^{2}\theta_{\min}}{(1+b)\mathcal{A}+c_1\tilde{g}}\right) \omega_{k}
        \end{align*}
        which implies that $(f(\mathbf{x}_{0})-f(\mathbf{x}_{k+1}))\rightarrow\infty$ whenever $k\rightarrow\infty$ using \eqref{cgspf1}. This is a contradiction because the sequence $\{\mathbf{x}_{k}\}$ belongs to a compact set, and $f$ is continuous. Therefore, our assumption is wrong, and hence we have $\lim\limits_{k\to\infty}\lVert \mathbf{g}_k\rVert=0$.
    \end{proof}
\end{theorem}
}

\section{Numerical Results}

INTHOP is implemented in C++. All computational runs were performed on a PC with a 12th Gen Intel Core i7-12700 processor (12 cores, 2.10 GHz) and 32 GB RAM. We used the \texttt{Gaol} library \cite{goualard2005gaol} to perform interval arithmetic operations. We employed the \texttt{Eigen} library for linear algebra operations, including Cholesky decomposition and inverse computation \cite{eigenweb}. Symbolic function representation and the computation of gradients and Hessians have been implemented using the \texttt{GiNaC} library \cite{bauer2002introduction}, which allows for symbolic differentiation in C++.
We have used the following parameters while implementing the above algorithms: gradient tolerance $\epsilon_g = 0.001$, maximum number of iterations $iter_{\max} = 10000$, backtracking factor $\nu = 0.5$, Armijo condition constant $\eta = 0.001$, initial Armijo step length $\theta_0 = 1$, regularization parameter for INTHOP:A1 $\beta = 1$, initial interval size $\delta_0 = 0.1$ and the interval scaling factor $r = 2$.  Parameters $c_1$ and $\tilde{g}$, to make the Hessian matrix positive definite, are chosen equal to  $0.001$ and $\lVert \textbf{g}_{t}\rVert$, respectively.

In this section, we apply the INTHOP algorithm on an extensive set of problems as referenced in \cite{problemset} whose dimensions range from 1 to 1404 variables. The complete problem set and the initial guess for each problem are given in ESI. These experiments evaluate the performance of the INTHOP variants compared to other algorithms. Specifically, we assess the convergence properties and computational efficiency of our method across the problem set, highlighting the advantages and trade-offs compared to the traditional methods that compute gradient and Hessian at each iteration.

\subsection{Benchmarking Interval Hessian Frameworks}
We utilize data profiles~\cite{doi:10.1137/080724083} to systematically benchmark the effectiveness of the optimization algorithm and investigate the effects of various factors on its performance. These profiles facilitate structured comparisons by explicitly defining three fundamental elements: a set of benchmark optimization problems, denoted by $\mathcal{S}$; a collection of optimization algorithms, denoted by $\mathcal{\tilde{A}}$; and an established convergence criterion, $\mathcal{T}$.

A data profile quantifies the absolute performance of an optimization algorithm. For each pair consisting of a problem $s \in \mathcal{S}$ and an algorithm $a \in \mathcal{\tilde{A}}$, we define a performance metric $m_{s,a}$, which, for example, could be the number of Hessian evaluations required to achieve the convergence criterion. Formally, the data profile $d_a(\zeta)$ of an algorithm $a$ indicates the proportion of benchmark problems that algorithm $a$ solves within a specified computational threshold $\zeta$. Specifically, it is given by:
\begin{equation}
    d_a(\zeta) = \frac{1}{|\mathcal{S}|}\text{size}\{s\in\mathcal{S}:m_{s,a}\leq \zeta\},
    \label{eq:dataprofile}
\end{equation}
where $|\mathcal{S}|$ denotes the total number of benchmark problems considered. When $m_{s,a} = \infty$, the algorithm, $a$ fails to satisfy the convergence criterion for problem $s$. A problem is considered solved if and only if the final point returned by an algorithm is a stationary point, i.e., it satisfies the following conditions: $\|\mathbf{g}_k\| < \epsilon_g$, where $\epsilon_g>0$. Data profiles are particularly valuable when a fixed computational budget is available, guiding the selection of algorithms most capable of solving a problem within the allocated computational budget.

\subsubsection{Effect of Interval Size}
A critical parameter that impacts the performance of INTHOP:F is the interval size $\delta$. A larger interval size can reduce the computational effort required for repeated Hessian evaluations and matrix inversion, but might lead to conservative steps that can be seen from the relation given in \eqref{lemma:regularization_effect} and \eqref{pklessinvres}. Conversely, smaller intervals offer accuracy in capturing curvature details but incur higher computational costs due to repeated Hessian evaluations. Additionally, the accuracy of $\alpha$ depends significantly on the specific estimation method used. Our observations indicate that the EM and MK methods provide tighter $\alpha$ than the GGN method. However, the two methods require $\mathcal{O}(n^3)$ operations, while the Gerschgorin method requires $\mathcal{O}(n^2)$ operations. To systematically investigate these trade-offs, we evaluate nine algorithmic variants, combining three distinct values (0.1, 0.5, 1) of $\delta$ with the three $\alpha$ estimation methods discussed in Section~\ref{sec:Methods}.



\begin{figure}[H]
    \centering
        \includegraphics[width=\linewidth]{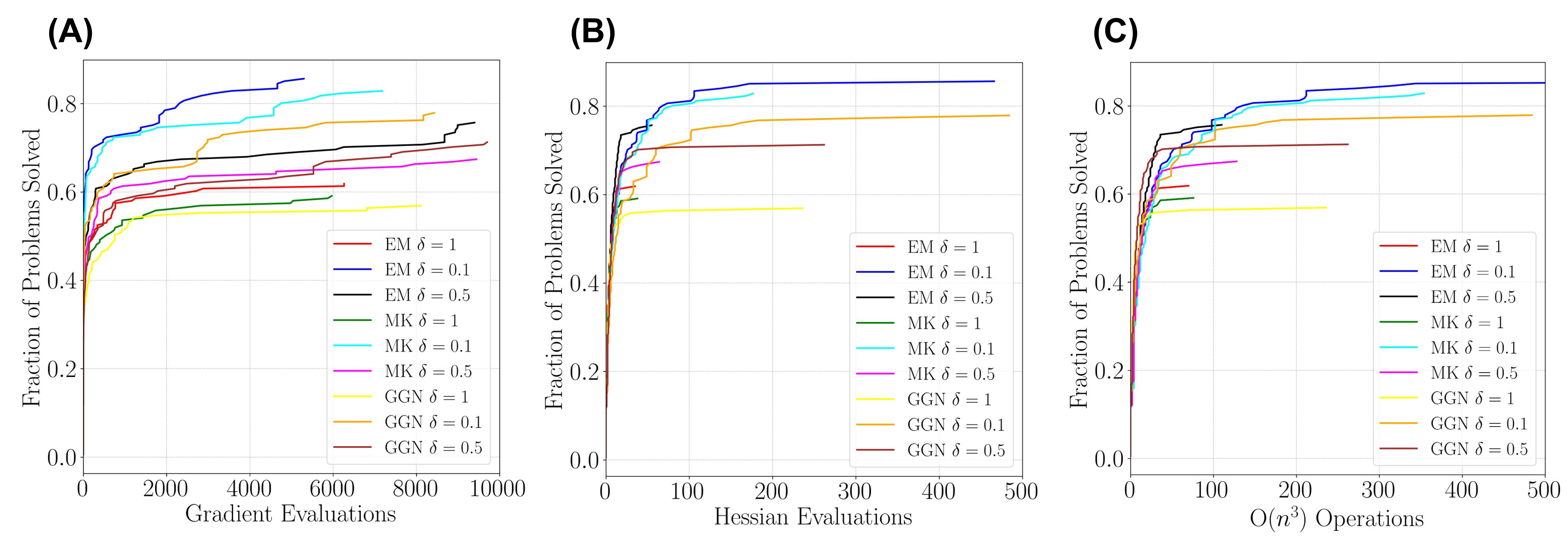}
    \caption{Fraction of problems solved by INTHOP:F variants with (A) gradient evaluations, (B) Hessian evaluations, and (C) $\mathcal{O}(n^3)$ operations} 
    \label{fig:Interval_Size_Comparison}
\end{figure}


These variants are tested across a comprehensive benchmark set comprising {181} optimization problems taken from \cite{problemset}. While the problem set consisted of bound-constrained problems, we solved them by considering them as unconstrained problems. The computational results are depicted in Figure~\ref{fig:Interval_Size_Comparison}, which shows the fraction of problems solved with the number of gradient and Hessian evaluations and $\mathcal{O}(n^3)$ operations. Each INTHOP variant is denoted in the figures using the nomenclature \texttt{<$\alpha$ calculation method>($\delta$=Interval size)}. Figure~\ref{fig:Interval_Size_Comparison} demonstrates that the proposed approach successfully solves approximately {80-85\% of the problems within 3000 gradient evaluations, 100 Hessian evaluations, and 200 $\mathcal{O}(n^3)$ operations}. Recall that if $\mathbf{x}_k\in [\mathbf{x}_t^L, \mathbf{x}_t^U]$, computing the search direction requires gradient evaluation and a matrix-vector product. In contrast, if $\mathbf{x}_k \notin [\mathbf{x}_t^L, \mathbf{x}_t^U]$, the Hessian and $\alpha$ are computed. The search direction is obtained by taking the product of $\nabla^2\mathcal{H}_t$ inverse with gradient at $\mathbf{x}_{k}$. Thus, contrary to the Newton method, which requires an equal number of gradient and Hessian evaluations to compute the search direction, INTHOP requires fewer Hessian evaluations than gradient evaluations. Furthermore, whenever $\mathbf{x}_k \notin [\mathbf{x}_t^L, \mathbf{x}_t^U]$, the MK and EM methods require two $\mathcal{O}(n^3)$ operations; one to estimate $\alpha$ and the second to compute the search direction. On the other hand, the GGN method requires one $\mathcal{O}(n^2)$ operation to compute $\alpha$ and one $\mathcal{O}(n^3)$ operation to compute the inverse of $\nabla^2\mathcal{H}_t$ to obtain the search direction. 

{We make the following five observations. First, INTHOP employing the EM method to estimate $\alpha$ performs the best, followed by MK and the Gerschgorin methods. Second, for the same $\alpha$ calculation method, increasing $\delta$ results in performance degradation across all metrics. Third, EM and MK methods with $\delta=0.1$ perform comparably across all performance metrics. Fourth, the EM method with $\delta=0.5$ performs comparably to the Gerschgorin method with $\delta=0.1$ for approximately 65\% of problems, which in turn achieves a performance superior to the MK method with $\delta=0.5$ in all performance metrics. Fifth, MK and Gerschgorin methods with  $\delta=0.5$ perform similarly in terms of gradient evaluations for approximately 65\% of problems.} However, the latter method outperforms in terms of Hessian evaluations and $\mathcal{O}(n^3)$ operations. The performance difference between these two methods is greater for $\mathcal{O}(n^3)$ operations, as the Gerschgorin method requires $\mathcal{O}(n^2)$ operations to compute $\alpha$, while the MK method requires $\mathcal{O}(n^3)$ operations. Our results illustrate that accurately estimating $\alpha$ is critical for the algorithm's performance, which depends on the calculation method of $\alpha$ and the size of the interval $\delta$.

{The INTHOP algorithm that employs EM method with $\delta=0.1$ could not solve 26 problems out of the set of 181 problems. For 5 problem instances, the algorithm reached the prescribed maximum iteration limit set at 10,000; for 18 problem instances, the step size became smaller than the tolerance; for 1 problem, the algorithm reached the maximum time limit; and for 2 problems, the algorithm could not find the direction after Hessian modification. The case when the algorithm is unable to find the descent direction is very rare. This behavior can be understood from the structure of the Hessian matrix. At some points, the Hessian becomes singular and does not provide curvature information in several directions. As a result, the system used to compute the search direction becomes difficult to solve or unstable. Even after modifying the Hessian, this issue may still remain, making it hard for the algorithm to find a reliable descent direction. }  

\vspace{5 pt}
\subsubsection{Adaptive Delta Evaluation}
In previous sections, $\delta$ was treated as a fixed hyperparameter and we observed that it significantly impacted the overall algorithmic performance. However, the optimal interval size is inherently problem-dependent and sensitive to local landscape variations. Larger intervals lead to a larger difference in the original Hessian and its approximation (Theorem \ref{theorem:theo1}). In comparison, smaller intervals improve tightness but incur higher computational costs due to frequent Hessian evaluations and $\mathcal{O}(n^3)$ operations.

\begin{figure}[H]
    \centering
    \includegraphics[width=\linewidth]{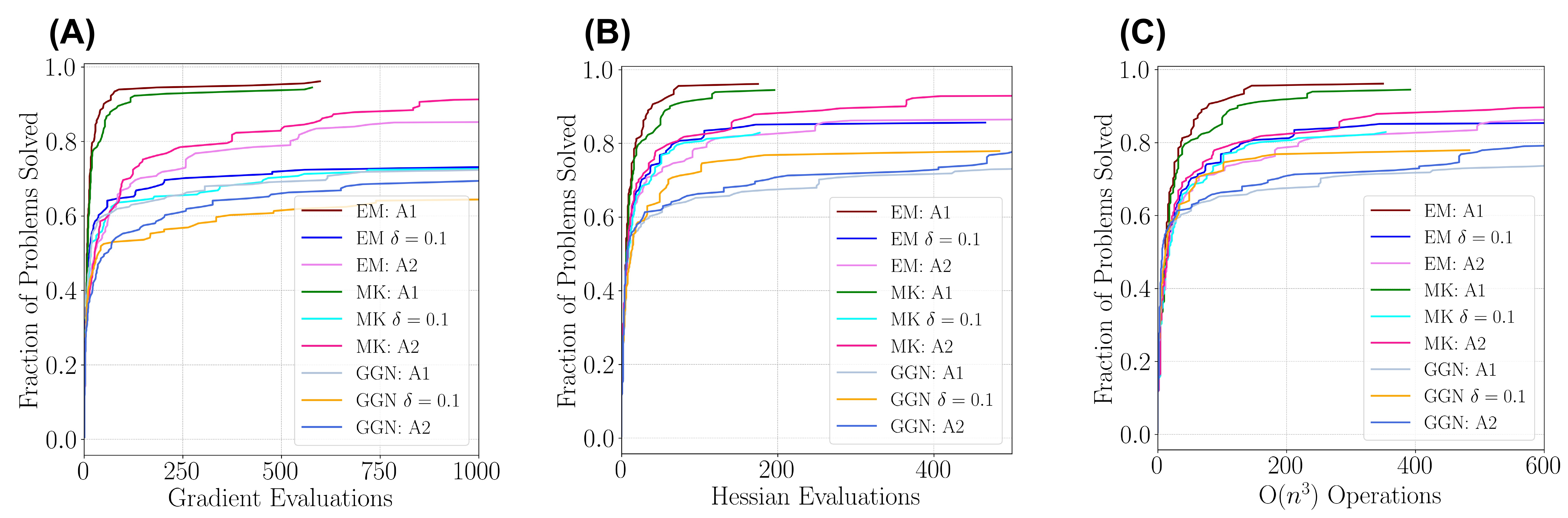}
    \caption{Data profile corresponding to fixed and adaptive $\delta$ INTHOP variants using (A) gradient evaluations, (B) Hessian evaluations, and (C) $\mathcal{O}(n^3)$ operations as the performance metrics}
    \label{fig:Adaptive_methods}
\end{figure}

To address this trade-off, we employ two adaptive interval strategies, INTHOP: A1 and INTHOP: A2, outlined in Algorithms~\ref{algo:AdaptGlobalStructure} and \ref{algo:TrustAdaptGlobalStructure}, respectively. These strategies dynamically adjust \(\delta\) to balance the Hessian approximation accuracy with computational overhead. The interval size is constrained within \([0.001, 10]\). The lower bound on $\delta$ ensures that the algorithm does not behave like a traditional point-based method, while the upper bound ensures that the Hessian approximations are not loose. To evaluate the effectiveness of these adaptive schemes, we combine INTHOP: A1 and INTHOP: A2 with each eigenvalue estimation method in Section~\ref{sec:Methods}, and compare them with the best performing INTHOP variant with fixed interval size corresponding to $\delta=0.1$, as illustrated in Figure \ref{fig:Interval_Size_Comparison}. Figure \ref{fig:Adaptive_methods} shows the resulting data profiles obtained for various methods. Each adaptive variant is denoted in the figures using the nomenclature \texttt{<Eigenvalue Method>:<Adaptive Variant>}. We make the following key observations. {First, for the same $\alpha$ calculation method, the INTHOP variant with adaptive methods to select $\delta$ performs superior to the best fixed $\delta$ method. Second, A1 with both methods, MK and EM, performs best and solve approximately 96\% of problems}. {{Third, EM: A1 was unable to solve 7 problems from the set of 181 problems. For 5 problem instances, the number of iterations exceeded, and for 2 problem instances, the algorithm could not find the direction.}} {Fourth, MK: A2 solves approximately 92\% of the problems, while EM: A2 solves  approximately 85\% of the problems. Fifth, for the GGN method, the A1 variant performs better than A2 with respect to gradient evaluations, but achieves a similar performance with respect to the number of Hessian evaluations and $\mathcal{O}(n^3)$ operations.}



\begin{figure}[H]
\centering
\includegraphics[width=\textwidth]{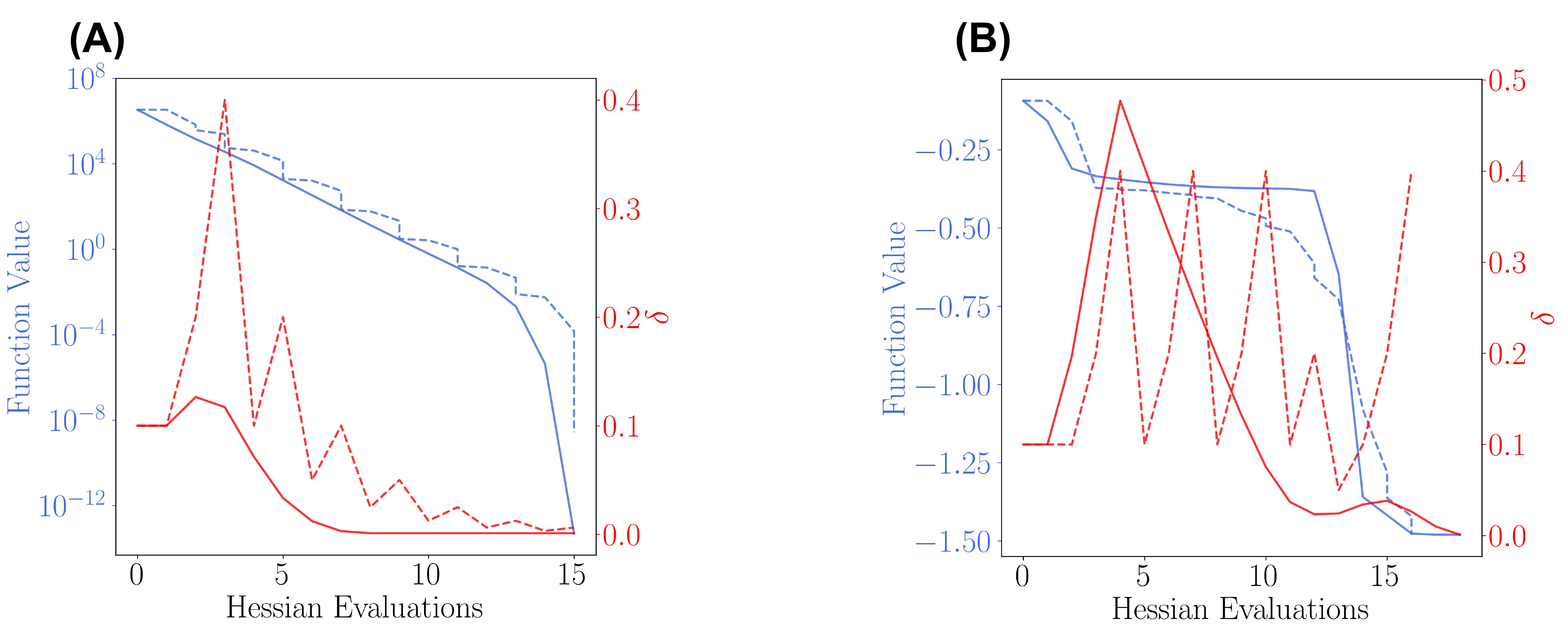}
\caption{Variation of function value and $\delta$ with iterations for (A) \texttt{6\_s273} and (B) \texttt{10\_Shekelfox10}.}
\label{fig:variation_delta}
\end{figure}

Another critical question is how adaptive variants modify $\delta$ values with iterations. To illustrate this, we show the variation of $\delta$ for two problems (\texttt{6\_s273} and {\texttt{10\_Shekelfox10}}) that were not solved by the fixed $\delta$ variant using the EM method to compute $\alpha$, while the two adaptive methods solved the two problems. The results for \texttt{6\_s273} and {\texttt{10\_Shekelfox10}} are shown in Figure \ref{fig:variation_delta} (A) and (B), respectively. The left y-axis shows the function values in blue, and the right y-axis shows $\delta$ in red with Hessian evaluations. The solid lines represent the performance of the adaptive method A1, and the dashed lines represent the A2 method. Recall that when $\textbf{x}_t\notin [\textbf{x}_k^L, \textbf{x}_k^U]$, $\delta$ is modified and the Hessian matrix is evaluated. For \texttt{6\_s273}, we observe that initially $\delta$ increases for both adaptive variants to allow a larger step size during the initial phase of the algorithm. However, $\delta$ decreases to a small value as the algorithm approaches the solution. For {\texttt{10\_Shekelfox10}}, for EM: A1, $\delta$ first increases then decreases strictly, whereas no apparent trend is observed for EM: A2. Importantly, we note that adapting $\delta$ can achieve faster convergence.


Overall, these results highlight the robustness of the adaptive interval strategies, particularly when used with the EM and MK $\alpha$ estimation methods. It consistently outperforms the fixed-interval approaches across all the performance metrics. 

\subsection{Comparison with Competing Algorithms}   
\begin{figure}[H]
\centering
\includegraphics[width=\textwidth]{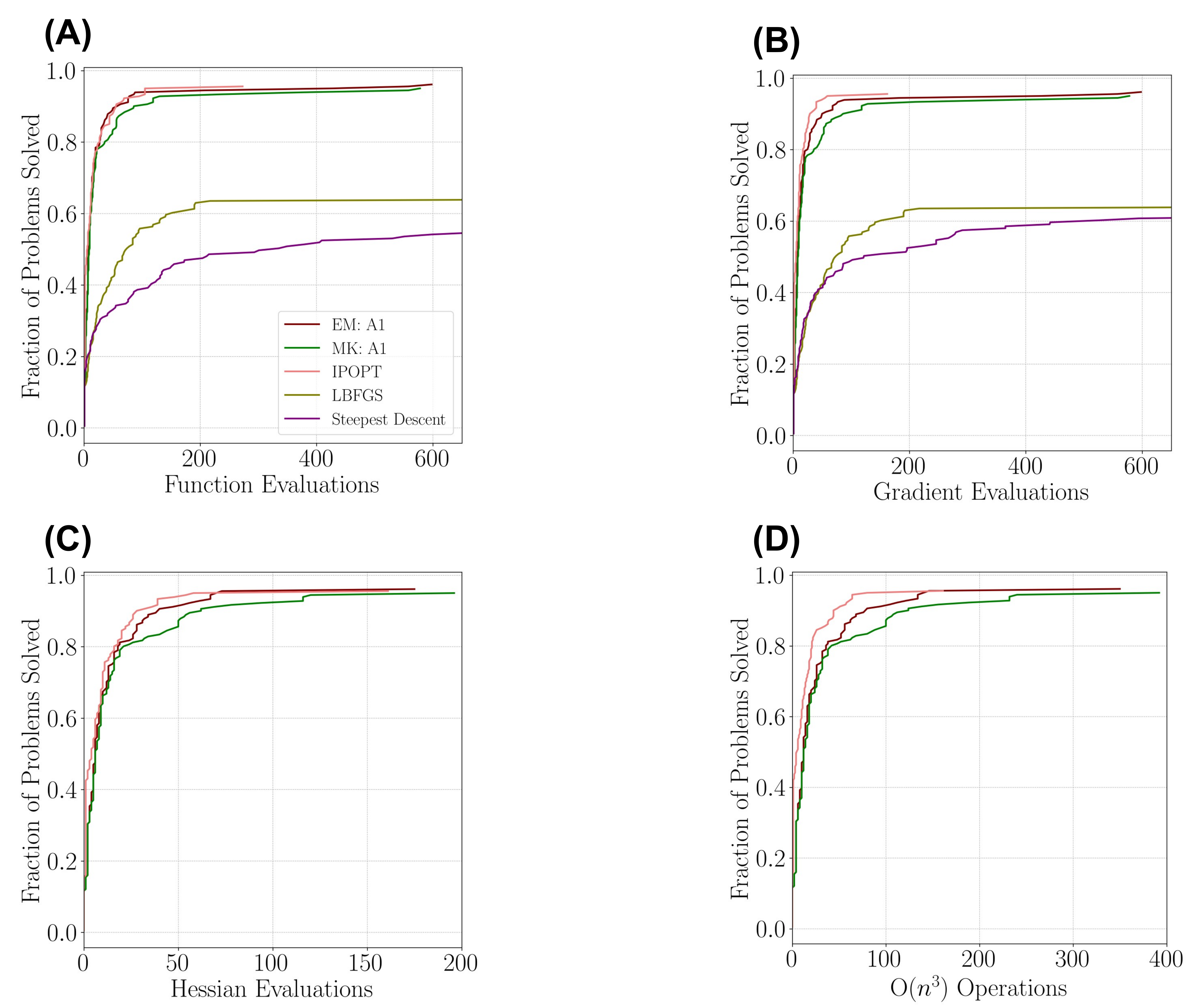}
    \caption{Fraction of problems solved by INTHOP with EM: A1 and MK: A1, steepest descent, LBFGS impletmented in PyTorch, and IPOPT solver with (A) function evaluations, (B) gradient evaluations, (C) Hessian evaluations, and (D) $\mathcal{O}(n^3)$ operations.}
    \label{fig:all-dp-comparison}
\end{figure}

{Next, we compare the performance of the proposed INTHOP algorithm with several well-established optimization methods, namely the IPOPT solver \cite{wachter2006implementation}, the steepest descent method \cite{nocedal1999numerical}, and the limited memory BFGS (LBFGS) algorithm implemented in PyTorch \cite{paszke2019pytorch}.

The data profiles are shown in Figures \ref{fig:all-dp-comparison} (A), (B), (C), and (D), showing the number of function evaluations, gradient evaluations, Hessian evaluations, and $\mathcal{O}(n^3)$ operations, respectively, required by the various methods. The steepest descent and LBFGS methods do not require Hessian evaluations and $\mathcal{O}(n^3)$ operations. Therefore, these two are not included in Figures (C) and (D). From these results, we draw the following key conclusions. 
First, the algorithms using Hessian information (EM: A1, MK: A1, and IPOPT) solve more problems in fewer function and gradient evaluations. On comparing the fraction of problems solved based on function evaluations, we observe that the INTHOP algorithms EM: A1 and MK: A1, and IPOPT solve approximately 80\% of problems within 50 function and gradient evaluations. In contrast, LBFGS and steepest descent methods underperform, struggling to solve many problems even with significantly higher function evaluations. Steepest descent, in particular, solves the fewest problems and requires the most function evaluations. This is due to two primary reasons: a larger number of iterations and repeated evaluations were needed to satisfy the Armijo condition for step length ($\theta$) selection. 

Third, INTHOP performs better than LBFGS, indicating that the approximate Hessian obtained by the former method is more accurate than the approximation of the Hessian in the latter. Fourth, IPOPT and EM: A1  solve the highest fraction of problems in all performance metrics. Although MK: A1 is not much behind these and solves almost the same number of problems. 

Furthermore, in terms of Hessian evaluations and $\mathcal{O}(n^3)$ operations, the proposed method remains highly competitive with IPOPT, indicating that it effectively balances second-order information with computational efficiency. These results clearly demonstrate that the proposed INTHOP algorithm provides a strong alternative to existing state-of-the-art methods, achieving high robustness and efficiency on a wide range of test problems.}



\section{Conclusion}
This article introduced INTHOP, an interval Hessian-based optimization framework for nonconvex optimization problems. The framework is based on approximating the original Hessian such that the approximate Hessian is positive definite in an interval. This strategy offers several advantages: (i) the search direction is guaranteed to be descent; (ii) the Hessian need not be recomputed for the iterates lying within the interval; (iii) the search direction can be obtained by matrix-vector product if the iterates are within the same interval. The key idea in finding a positive definite approximation of the Hessian matrix is to estimate the minimum eigenvalue of the interval Hessian matrix. The algorithm performance critically depends on the size of  the interval and the method employed to compute
the minimum eigenvalue of the interval Hessian. Accordingly, we develop INTHOP variants where the interval size remains the same, and another variant where the interval size changes with iterations. We prove the global convergence that is applicable to all the proposed variants of INTHOP.  Furthermore, we implement several methods of $\mathcal{O}(n^2)$ and $\mathcal{O}(n^3)$ complexities to compute the minimum eigenvalue of the interval Hessian matrix. We compare the performance of INTHOP variants on a set of benchmark problems and demonstrate that INTHOP with adaptive interval size and using an $\mathcal{O}(n^3)$ complexity method outperforms the other INTHOP variants across multiple metrics, including function, gradient, and Hessian evaluations and $\mathcal{O}(n^3)$ operations. 

{We also compare the performance of INTHOP with classical point-based methods, including the steepest descent method, the L-BFGS algorithm implemented in PyTorch, and the state-of-the-art IPOPT solver, across a diverse set of benchmark problems, including large-scale instances with high dimensional variables. The results indicate that INTHOP requires significantly fewer function and gradient evaluations than the steepest descent and L-BFGS. In particular, the steepest descent and LBFGS exhibit slower convergence and require a larger number of iterations due to repeated function and gradient evaluations for the line search procedures.

In comparison with IPOPT, which effectively exploits second-order information, INTHOP achieves competitive performance across all metrics (function evaluation, gradient evaluation, Hessian evaluation, and number of $\mathcal{O}(n^{3})$ operations).

Moreover, it has also been observed that for some highly nonconvex problems, INTHOP often requires substantially fewer $\mathcal{O}(n^3)$ operations than IPOPT, highlighting its computational advantage. Overall, these results highlight that the proposed INTHOP algorithm offers a robust, efficient, and competitive alternative to existing state-of-the-art optimization methods across a diverse range of test problems.}

A notable weakness of INTHOP is that it does not exploit the exact Hessian in convex regions of a nonconvex function, potentially missing performance gains in some problem instances. We aim to address this limitation in our future work and extend INTHOP to solve constrained optimization problems.

    \section*{Data availability}

The data is included as electronic supplementary information.

\section*{Conflict of interest}

 The authors declare that they have no known competing financial interests or personal relationships
that could have appeared to influence the work reported in this paper.

\addcontentsline{toc}{section}{References}
\bibliographystyle{apalike}
\bibliography{ref}


\clearpage 
\appendix
\addcontentsline{toc}{section}{Appendix} 
\appendix

\section{Proofs} \label{appendix:proofs}

This appendix provides proofs of Theorem \ref{theorem:theo1} and Lemma \ref{lemma:regularization_effect} stated in Section  \ref{sec:SearchDirectionConstruction}. 
We prove the following results stated in Theorem \ref{theorem:theo1}:

\begin{itemize}
    \item \textbf{Function value bound:}
       \begin{align*}
    \left| \mathcal{L}(\mathbf{x}_t) - f(\mathbf{x}_k) \right| \leq \frac{L_f}{2} \sqrt{n} \delta + \frac{|\lambda_{min}|}{8} n \delta^2
\end{align*}
    
    \item \textbf{Gradient bound:}
    \[
    \left\| \nabla \mathcal{L}(\mathbf{x}_t) - \nabla f(\mathbf{x}_k) \right\| \leq \frac{L_g}{2} \sqrt{n} \delta + \frac{|\lambda_{min}|}{2} \sqrt{n}\delta
    \]
    
    \item \textbf{Hessian bound:}
    \[
    \left\| \nabla^2 \mathcal{L}(\mathbf{x}_t) - \nabla^2 f(\mathbf{x}_k) \right\| \leq \frac{L_H}{2} \sqrt{n} \delta + \frac{|\lambda_{\min}|}{2} \sqrt{n} 
    \]
\end{itemize}

Proofs of each of the above results are provided below.

\noindent \textbf{Proof of function value bound:}
The difference between the underestimator ($\mathcal{L}$) and the true function ($f$) satisfies
\begin{align*}
    \left| \mathcal{L}(\mathbf{x}_t) - f(\mathbf{x}_k) \right| 
    &= \left| f(\mathbf{x}_t) - f(\mathbf{x}_k) + \sum_{i=1}^n \alpha (x_i^L - x_{i,t})(x_i^U - x_{i,t}) \right|.
    \end{align*}
Using triangle inequality, we get 
    \begin{align*}
       \left| \mathcal{L_f}(\mathbf{x}_t) - f(\mathbf{x}_k) \right| 
       &\leq \left| f(\mathbf{x}_t) - f(\mathbf{x}_k) \right| + \sum_{i=1}^n \alpha \, |(x_i^L - x_{i,t}) (x_i^U - x_{i,t})|.
       \end{align*}
Using Eq. \ref{eq:lipsch_f}, the fact that the maximum of $|(x_i^L-x_{i,t})(x_i^U-x_{i,t})|$ is at $x_{i,t}=\frac{x_i^U+x_i^L}{2}$ and Eq. \ref{eq:alpha1}, we get 
\begin{align*}
    &\leq L \|\mathbf{x}_t - \mathbf{x}_k\| + \sum_{i=1}^n \frac{|\lambda_{min}|}{8} (x_i^U-x_i^L)^2 .
\end{align*}
We note that $||\mathbf{x}_t-\mathbf{x}_k|| \leq \sqrt{n}\frac{\delta}{2}$ and $x_i^U-x_i^L=\delta$. Using in the above equation, we get the desired result.
\qed

\vspace{1em}

\noindent \textbf{Proof of gradient bound:}
Similarly, for the gradient difference, we have
\begin{align*}
    \left\| \nabla \mathcal{L}(\mathbf{x}_t) - \nabla f(\mathbf{x}_k) \right\| 
    &= || \nabla f(\mathbf{x}_t) - \nabla f(\mathbf{x}_k) + \alpha [2\mathbf{x}_{t}-(\mathbf{x}^L+\mathbf{x}^U)]|| \\
    &\leq \left\| \nabla f(\mathbf{x}_t) - \nabla f(\mathbf{x}_k) \right\| +  |\alpha| \, \|\mathbf{x}^L - \mathbf{x}^U\| \\
    &\leq L_g \|\mathbf{x}_t - \mathbf{x}_k\| + |\alpha| \sqrt{n} \delta \\
    &\leq \frac{L_g}{2} \sqrt{n} \delta + \frac{|\lambda_{min}|}{2} \sqrt{n}\delta.
\end{align*}
\qed

\vspace{1em}
\noindent \textbf{Proof of Hessian Bound:}
Finally, for the Hessian difference, we have
\begin{align*}
    \left\| \nabla^2 \mathcal{L}(\mathbf{x}_t) - \nabla^2 f(\mathbf{x}_k) \right\| 
    &= \left\| \nabla^2 f(\mathbf{x}_t) - \nabla^2 f(\mathbf{x}_k) + \alpha \mathbf{I}\right\| \\
    &\leq \left\| \nabla^2 f(\mathbf{x}_t) - \nabla^2 f(\mathbf{x}_k) \right\| +  |\alpha| \sqrt{n} \\
    &\leq L_H \|\mathbf{x}_t - \mathbf{x}_k\| + |\alpha| \sqrt{n} \\
    &\leq \frac{L_H}{2} \sqrt{n} \delta + \frac{|\lambda_{\min}|}{2} \sqrt{n}.
\end{align*}
\qed

\end{document}